\documentclass{amsart}
\usepackage{amsfonts,amsthm,amsmath}
\usepackage{amssymb}
\usepackage{amscd}
\usepackage[utf8]{inputenc}
\usepackage[a4paper]{geometry}
\usepackage{enumerate}
\usepackage[usenames,dvipsnames]{color}
\usepackage{array}
\usepackage{array,tabularx}

\usepackage{amsthm}
\usepackage{tikz}
\usepackage{xcolor}

\DeclareMathOperator{\Coeff}{Coeff}
\DeclareMathOperator*{\Res}{Res}
\newcommand\note[1]{\mbox{}\marginpar{ \scriptsize\raggedright
\hspace{1pt}\color{red} #1}}

\numberwithin{equation}{section}
\numberwithin{equation}{subsection}

\theoremstyle{plain}

\newtheorem{theorem}[equation]{Theorem}
\newtheorem{lemma}[equation]{Lemma}

\newtheorem{corollary}[equation]{Corollary}

\newtheorem{thm}[equation]{Theorem}

\theoremstyle{definition}
\newtheorem{example}[equation]{Example}
\newtheorem{remark}[equation]{Remark}

\newtheorem{definition}[equation]{Definition}

\newtheorem{bekezdes}[equation]{}

\def\Z{\mathbb Z}

\newcommand{\calv}{{\mathcal V}}

\newcommand{\calt}{{\mathcal T}}

\newcommand{\calO}{{\mathcal O}}

\newcommand{\calS}{{\mathcal S}}

\newcommand{\calL}{\mathcal{L}}

\newcommand{\tX}{\widetilde{X}}

\newcommand{\cX}{{\mathcal X}}
\newcommand{\cO}{{\mathcal O}}
\newcommand{\bP}{{\mathbb P}}
\newcommand*{\linebundle}{\mathcal{L}}
\newcommand{\bC}{{\mathbb C}}
\newcommand{\cF}{{\mathcal F}}

\newcommand{\pic}{{\rm Pic}}

\newcommand{\bt}{{\mathbf t}}

\newcommand{\bZ}{{\mathbb{Z}}}
\newcommand{\bQ}{{\mathbb{Q}}}

\author{J\'anos Nagy and Andr\'as N\'emethi}
\address{HUN-REN Alfr\'ed R\'enyi Institute of Mathematics, Re\'altanoda utca 13--15., H-1053 Budapest,  Hungary}
\email{nagy.janos@renyi.hu}
\email{nemethi.andras@renyi.hu}

\thanks{NJ was partially supported by the KnotSurf4d ERC Advanced Grant of A. Stipsicz;
NA was partially supported by Hungarian National Research, Development and Innovation Office NKFIH Grant ``\'Elvonal (Frontier)'' KKP 144148. }

\title{Lipschitz inner rates of normal surface singularities}

\begin{document}

\keywords{normal surface singularities, inner rate, Lipschitz invariants, analytic semigroup}

\subjclass[2010]{Primary. 32S05, 32S25, 32S50, 57M27
Secondary. 14Bxx, 14J80, 57R57}

\begin{abstract}
In this article, we connect the local Lipschitz geometry of a normal surface singularity with a rational homology sphere link to the analytic semigroup $S_{an}(\tilde{X})$ associated with a convenient  resolution. Under certain restriction,  we provide an explicit computational characterization of the Lipschitz rate $q_v$
associated with an irreducible exceptional curve $E_v$ in terms of the semigroup. We also show that removing  the restriction  we get examples when $q_v$ cannot be determined from the semigroup. Interestingly  enough, this phenomenon already happens for rational  singularities. Indeed, we exemplify that the inner rates of rational singularities are not topological.
\end{abstract}

\maketitle

\linespread{1.2}


\pagestyle{myheadings} \markboth{{\normalsize  J. Nagy and A. N\'emethi}}
 {{\normalsize Lipschitz invariants}}

\section{Introduction}
\label{s:intro}

\subsection{}
The study of the local Lipschitz geometry of singular spaces has its roots in the Zariski saturation theory and the pioneering work of Pham and Teissier \cite{pham-teissier} on complex curve germs. Later, this framework was significantly expanded by the foundational contributions of several authors, see e.g. the work of 
Mostowski \cite{Mostowski},  Parusi\'nski \cite{parusinski} or Risler and Trotman \cite{Risler}  regarding Lipschitz stratifications of subanalytic sets and the metric behavior of singular spaces.
Then the  theory has undergone rapid development 
with several outstanding  results and connections, see the work of Birbrair, Ruas, L\^e, Fernandes, Sampaio, Henry, Parusi\'nski,  de Bobadilla, Neumann, Pichon,  and others in the subject --- for a list a detailed references see the  review article of  Ruas  \cite{Ruas} or of  Pichon 
\cite{Pichon}, or the other articles of the LNM volume \cite{springer-book}  and the references therein.

 For complex analytic surface singularities, this field has witnessed a remarkable activity and  fulfillment  in recent years. Initiated by the seminal works of Birbrair, Fernandes, Neumann, and Pichon \cite{bfn2008, bnp2014}, the classification of singularities up to inner and outer Lipschitz equivalence has emerged as a major theme connecting algebraic geometry, topology, and metric geometry of normal surface singularities. A key milestone in this development was the introduction of the thick-thin decomposition by Birbrair, Neumann, and Pichon \cite{bnp2014}, which provided a complete geometric description of the inner  metric structure.
More recently,  de Bobadilla, Heinze, Pe Pereira, and  Sampaio introduced a new metric homology theory, called the 
 `moderately discontinuous homology', which captures Lipschitz properties of metric subanalytic germs \cite{bob1}, see also \cite{bob2}. 
All these works establish deep algebraic bridge connecting  metric invariants vith specialized differential complexes on singular varieties.

 Within this broad context, a crucial tool in analyzing these structures remains the concept of 
{\bf contact exponents (or rates) between curve germs, which measures how fast two curves diverge} with respect to either the inner metric $d_i$ (induced by lengths of paths contained within the singularity) or the outer metric $d_o$ (induced by distance in the ambient Euclidean space).

\subsection{}
Let $(X, 0)$ be a germ of a complex analytic normal surface singularity and $\pi : \tilde{X} \to X$ a good resolution with exceptional curve $E = \cup_{v \in \mathcal{V}} E_v$. 
 The topological type of $(X,0)$  is characterized by the  link $\partial \tilde{X}$ (or by the dual resolution graph $\mathcal{T}$), while the analytic type by the  local
 algebra $\mathcal{O}_{(X,0)}$. 

A crucial bridge between the metric geometry   and the analytic structure of $(X, 0)$ is provided by generic transversal curvettes $C_1, C_2 \subset \tilde{X}$ meeting a
fixed  irreducible exceptional component $E_v$ at distinct generic points,  and  by the contact properties of their projections 
$D_i=\pi(C_i)\subset X$, $i=1,2$.
  By a recent theorem of Belotto da Silva, Fantini and Pichon \cite{fantini}, the inner and outer contact exponents $q_i^X(D_1, D_2)$ and $q_o(D_1, D_2)$ of the pair  $D_1, D_2 \subset X$ coincide and depend solely on the vertex $v$, defining a well-defined metric invariant,  denoted by $q_v$, and called `inner rate of $v$'. 

{\bf The set of inner rates $\{q_v\}_v$ play crucial role in  Lipschitz geometry  classification}. E.g.,  the inner Lipschitz geometry of $(X,0)$ is uniquely determined by the set of 
inner rates $\{q_v\}_v$  and the maximal ideal cycle $Z_{max}$
 \cite{bnp2014}. 
Here $Z_{max}$ is that part of the divisor of the generic linear function which is supported on $E$. 

\subsection{} 
{\bf The primary goal of this article is to provide an explicit computational characterization of the Lipschitz invariants $\{q_v\}_v$ in terms of the analytic semigroup }
$\mathcal{S}_{an}(\tilde{X}) = \{\operatorname{div}_E(\pi^*f) \mid f \in \mathcal{O}_{(X,0)}\}$ associated with a convenient resolution $\tX$. 
In this way we create a new bridge between the Lipschitz geometrical 
invariants $\{q_v\}_v$ and the  discrete analytic invariant $\mathcal{S}_{an}(\tX)$.
By establishing such  connection, we show how the analytic structure of the semigroup $
\mathcal{S}_{an}(\tilde{X})$ directly governs the metric rates.

In our discussions we will assume that {\it the link  is a rational homology sphere}. This is equivalent with the fact that 
 the dual resolution graph $\mathcal{T}$  is a tree and $E_v$'s  are rational curves.

\subsection{}
{\bf It is natural to  investigate whether $q_v$ is a purely topological invariant} determined by the resolution graph $\mathcal{T}$. 
In other words, we ask whether the inner rates are the same for all analytic structures supported by a fixed  topological type. 
Or, we can ask its weaker, local version,  whether $q_v$  stays constant under a local analytic  topologically constant deformation of surface singularities.

The answer for both cases is negative, specific examples are known  for both cases. 
If one accepts  non rational homology sphere links too (as the literature usually does), then 
one can construct easily a topological type which supports  two different analytic structures, which have
different geometric genera, different analytic semigroups, and different inner rates. Regarding the topological constant deformations, along  the  
 classical Brian\c{c}on--Speder  deformation $x^5 + z^{15} + y^7x + t y^6 z^3 = 0$ \cite{BS} the link is constant, but  some of the Lipschitz invariants $q_v$ change.
In fact,  this family is not Zariski equisingular, which suggest that the Lipschitz geometry also must jump. Note that the constant link in the 
 Brian\c{c}on--Speder  deformation is not a rational homology sphere. 

Since in this article we treat the case of singularities with  rational homology sphere links, using our results we can exemplify the above two negative statements for such singularities as well. We provide an integral homology sphere link, and two analytic types supported by it, with different $q_v$ (for a certain  well-chosen $v$).
Also, we provide a  Brian\c{c}on--Speder  type deformation, with constant rational homology speher link, such that along the deformation both the maximal ideal cycle and some of the inner rates jump. (In both cases the analytic semigrours jump as well.) For details see subsection \ref{ss:examples}. 
 However, these examples do not refute the possible connection between the analytic semigroup and the internal rates.

Having such counterexamples in the general cases, we can also  try to restrict ourselves to specific families of singularities. For example, in the case of 
rational singularities several analytic invariant are topological, e.g. the maximal ideal cycle, the analytic semigroup, cohomology groups of line bundles on the resolution,
see e.g. \cite{Artin62,Artin66, Lauferbook,Lipman,NBook}.  
 Hence, it is natural to ask: are the inner rates also topological?
{\bf This question was our second objective: to prove the topological nature of inner rates for rational singularities.}
 As we will see below, this expectation does not always come true!

\subsection{The main results} 
In the body of the paper,  we will show that under some assumptions,  for general (not necessarily rational)  singularities,
 the inner rate can indeed be computed from the analytic semigroup. In particular, for rational singularities, when these assumptions hold,  the inner rates are indeed topological.
However, in general both statement fail. The geometrical picture behind the validity or the failure of the statements is the following.
When the multiplicity of $Z_{max}$ along $E_v$ is $m>1$ then {\bf $q_v$ is strongly related with the properties of the 
$m$-fold covering  of $(X,0)$ branched along the zero set of the generic linear function}. The point is that the geometry of this new singularity $(X^{(m)},0)$ is more complicated than
the geometry of $(X,0)$. When the needed information formulated at the level of $(X^{(m)},0)$ can be reinterpreted at the level of $(X,0)$ then we have positive results. Otherwise we obtain counterexamples. Surprisingly, even if $(X,0)$ is rational, the analytic semigruop of $(X^{(m)},0)$ can not be deduced from the topology of $(X,0)$.
We can construct a rational  topological type, with two different analytic types,  which provide for a certain $E_v$ two different inner rates.

In the main construction,  one of the irreducible exceptional  curves  of $\tX$ is covered by an elliptic curve in the resolution of $(X^{(m)},0)$, with an interesting 
intersection pattern with its neighbors. Varying these intersection points (in some linear system) 
we keep the topology fixed, however we can modify the analytic type in such a way that 
$q_v$ jumps.  The picture is really transparent at the level of $(X^{(m)},0)$ . 

This proves that the  Lipschitz geometry remains a fine analytic invariant even within the category of rational resolution graphs, showing clear sensitivity to the
analytic moduli of the plumbing construction.

\subsection{More detailed presentation of the main results} 
 We assume all the time that 
the resolution  $\tX$ factors through the blow up of the maximal ideal.

 It is known (see e.g.  \cite{fantini}, but it is reproved here as well)
 that in such a case, for any fixed $E_v$, if $(Z_{max}, E_v)<0$ then $q_v=1$. (Here $(\ ,\ )$ is the intersection form on $\tX$.)
 Hence  the interesting geometry lies in the complementary case $(Z_{max}, E_v) = 0$. 

In the sequel we fix a vertex $v$ of $\mathcal{T}$.
Below, in cases (1) and (3),  when we provide explicit expression for $q_v$ in terms of the analytic semigroup of $\tX$,  we start with a resolution $\tX$ which 
 factors through the blow up of the maximal ideal, and we also modify it by blowing up the intersection points of $E_v$ with its neighbors sufficiently many times
(the number of needed blow ups is bounded explicitly, see subsection \ref{ss:AS} and \ref{bek:modif}). 
We will continue to denote this modified  resolution  by $\tX$.

{\bf Our new main contributions in the case $(Z_{max}, E_v)=0$  can be summarized as follows:}

\begin{enumerate}
  \item \textbf{The $E_v$-multiplicity of $Z_{max}$, $Z_{max,v}$, equals one:} In Section 3, we establish several expressions and 
 a purely combinatorial characterization of $q_v$ from the analytic semigroup $\mathcal{S}_{an}(\tX)$ whenever  $(Z_{max}, E_v) = 0$ and $Z_{max,v} = 1$. 
For  concrete expression see Theorem \ref{innerrateone} and Corollary \ref{cor:case1}.  

    It is worth noting that this multiplicity one setup provides a direct generalization of the classical minimal rational case. Indeed, our general formula successfully recovers the explicit computation of the rates $q_v$ obtained by Birbrair, Neumann, and Pichon \cite{minimal-lne} for minimal rational surface singularities, where the condition $Z_{max,v} = 1$ holds 
automatically for every $v$ from the very definition of the family of minimal rational germs. 

\smallskip

    \item \textbf{The general case of  $m:=Z_{max,v} > 1$ and  cyclic coverings:} To tackle the case of $(Z_{max}, E_v)=0$ and  $m := Z_{\max,v} > 1$, we use the technique of  ramified cyclic coverings.
 We consider  the cyclic $\Z_m$-covering $(X^{(m)},0) \to (X,0)$ ramified along the zero set  of the generic linear section. We also consider the normalized fiber product 
$\tX^{(m)}$ of the covering $(X^{(m)},0) \to (X,0)$  with the resolution $\tX\to X$, and we consider the corresponding  covering $\tX^{(m)}\to \tX$ as well. 
Let $\{E^{(m)}_{v_i}\}_i $ be the  irreducible exceptional curves  of $\tX^{(m)}$ which cover $E_v$, they are permuted by the cyclic action, and their 
inner rates  $q^{(m)}_{v_i}$ (as invariant of  $(X^{(m)},0)$) are independent of the choice of $v_i$. 
We show that the inner rate $q_{v_i}^{(m)}$ (associated with $(X^{(m)},0)$ and $v_i$) and the original inner rate $q_v$ satisfy the scaling relation $q_{v_i}^{(m)} = m \cdot q_v$.
See Theorem \ref{ciclicinner}. 

It turns out that the maximal ideal cycle $Z^{(m)}_{max}$ of $\tX^{(m)}$ has multiplicity one along  $E^{(m)}_{v_i}$. Hence, in this way, 
for the singularity $(X^{(m)},0)$  and its resolution $\tX^{(m)}\to X^{(m)}$ we can apply part (1). Hence,  $q_{v_i}^{(m)} $, in 
particular, $q_v$ too, can be recovered from the 
analytic semigroup of  $\tX^{(m)}$. 
(However, this fact does not necessarily imply that they can be recovered from the semigroup of $\tX$ too!)


\smallskip

    \item \textbf{$m:=Z_{max,v} > 1$,  the `unbranched' case:}
Assume that we are in the situation of (2) but we assume additionally that above $E_v$ there is only one exceptional curve in $\tX^{(m)}$
(that is,  $\{E^{(m)}_{v_i}\}_i $ has cardinality one). Then,  $q_{v_1}^{(m)} $  considered in part (2) --- hence $q_v$ too, cf. part (2) ---
  can be computed in terms of the semigroup of $\tX$ as well.
For the explicit formulae see Thorems \ref{th:4.2.6} and \ref{th:4.2.7}.   


\smallskip

    \item  \textbf{$m:=Z_{max,v} > 1$,  the `branched' case:}  Finally, we show that in the general case of $|\{E^{(m)}_{v_i}\}_i|>1 $,
the inner rate cannot be computed from the semigroup $\mathcal{S}_{an}(\tX)$. 
In Remark \ref{rem:key} we highlight the point that distinguishes the  $|\{E^{(m)}_{v_i}\}_i|=1 $ case from the general one. What makes the discussion even more special and 
remarkable is the fact that this new phenomenon can appear already at the level of  rational singularities. 
 
Indeed, in Section 5  we prove the 
{\bf `analytic dependency of $q_v$ in the case of  rational singularities'}: 
We  provide  a family of analytic structures,  supported on a fixed  rational topological type,  such that along the deformation the 
inner rate associated with a well--chosen vertex jumps. Since for rational singularities the analytic semigroup is topological, this case shows that, 
in general, for arbitrary singularities,  the inner rate cannot be computed from the analytic semigroup  --- and this statement, rather surprisingly,  is valid already in the rational 
case. In fact,  if for an arbitrary singularity, we wish to compute successfully $q_v$ in terms of a certain analytic invariant,  that invariant must be non-topological for rational singularities. Of course, rational singularities have in general non-trivial moduli, however, to find such analytic invariants is really a challenge. 


This shows the subtlety of the invariant $q_v$ rather well. 
\end{enumerate}


\section{Preliminaries}\label{s:prel}

\subsection{The resolution}\label{ss:notation}
Let $(X,0)$ be a germ of a complex analytic normal surface singularity and  $\pi:\widetilde{X}\to X$  a good resolution of $(X,0)$. We denote the exceptional curve $\pi^{-1}(0)$ by $E$, and let $\cup_{v\in\calv}E_v$ be its decomposition into irreducible components.


Let $\calt$ be the dual resolution graph associated with $\pi$; it is a connected graph. Then $M:=\partial \widetilde{X}$ can be identified with the link of $(X,0)$, it is an oriented plumbed $3$--manifold associated with $\calt$. {\it We will assume that $M$ is a rational homology sphere,} or, equivalently, $\mathcal{T}$ is a tree and all genus decorations of $\mathcal{T}$ are zero. We use the same notation $\mathcal{V}$ for the set of vertices of $\calt$.

The homology group $L:=H_2(\widetilde{X},\mathbb{Z})$, endowed with a negative definite intersection form $(\,,\,)$, is a lattice. It is freely generated by the classes of oriented $2$--spheres $\{E_v\}_{v\in\mathcal{V}}$.
%
%
We define the {\it Lipman cone} (or, the topological semigroup)  as $\calS:=\{l\in L\,:\, (l, E_v)\leq 0 \ \mbox{for all $v$}\}$.

The lattice $L$ is partially ordered: $l'\geq l''$ whenever $l'_v\geq l''_v$ for all $E_v$-coordinates.

For any $l\in L$ there exists a unique minimal cycle $s(l)$   with $s(l)\in \calS$ and $s(l)\geq l$. It can be found by the following  Generalized Laufer Algorithm \cite[Prop. 6.6.3]{NBook}
which
constructs inductively a computation sequence $\{x_i\}_{i=0}^t$, $x_i\in L$, $x_{i+1}=x_i+E_{v(i)}$ for some $v(i)\in\mathcal{V}$.   Set $x_0:=l$. Assume that $x_i$ is already constructed. If $(x_i,E_v)\leq 0$ for every  $v$, then stop; $s(l)$ is this $x_i$.
 Otherwise, set $x_{i+1}=x_i+E_v$ where $E_v$ is (a choice) with
 $(x_i,E_v)>0$, and continue the test with $x_{i+1}$.

For any $f\in \calO_{(X,o)}$ let ${\rm div}(\pi^*f)$ be the divisor of $\pi^*f$ in $\tX$, it can be written as ${\rm div}_E(\pi^*f)+{\rm St}_{\tX}(f)$, where
${\rm div}_E(\pi^*f)$ is supported on $E$ and ${\rm St}_{\tX}(f)$ is the strict transform of $\{f=0\}$.  We define the {\it analytic semigroup} of $(X,o)$ associated with $\pi$ as
 $\mathcal{S}_{an}:=\{{\rm div}_E(\pi^*f)\, | \, f\in \calO_{(X,0)}\}$. Clearly, $\mathcal{S}_{an}\subset \mathcal{S}$.

\begin{lemma}\label{lem:iv} For any fixed resolution  $\pi$ and vertex $v$ of $\mathcal{T}$
there exists $r_v\in \mathbb{Q}_{>0}$ (depending only on $\mathcal{T}$ and $v$)
 such that
$$({\rm St}_{\tX}(f) , E_v) \leq r_v\cdot {\rm div}_E(\pi^*f)_v \ \ \mbox{for any $f\in \calO_{(X,0)}$}.$$
\end{lemma}
\begin{proof} Since ${\rm div}_E(\pi^*f)+{\rm St}_{\tX}(f)=  {\rm St}_{\tX}(f)+\sum_u  {\rm div}_E(\pi^*f)_u E_u$
is numerically trivial, 
$0=( {\rm St}_{\tX}(f), E_v) +   {\rm div}_E(\pi^*f)_v ( E_v, E_v) +\sum_{u\not=v}  {\rm div}_E(\pi^*f)_u( E_u, E_v)\geq 
( {\rm St}_{\tX}(f), E_v) +   {\rm div}_E(\pi^*f)_v ( E_v, E_v)$. Hence $r_v=-E_v^2$ is a possible value
(however, there are several other  choices as well). 
%
\end{proof}

The {\it maximal ideal cycle} $Z_{max}\in L$ of $\pi$ is defined as ${\rm div}_E(\pi^*L)$, where $L$ is a generic linear function in the   maximal ideal
$\mathfrak{m}_{(X,0)}$  of $\mathcal{O}_{(X,o)}$.

\subsection{Outer and inner Lipschitz invariants of a normal surface singularity}\label{ss:1.7}\cite{fantini,lne-surfaces,bnp2014,neumpichon-char,neumpichoncurve}

We consider  the following equivalence relation:
given two function germs $f, g \colon ([0, \infty), 0) \to ([0, \infty), 0)$,
 we say that  $f$ is equivalent with  $g$, and we write $f(t) = \Theta(g(t))$, if there exist real numbers $\eta > 0$ and $K > 0$ such that for all
  $t \leq \eta$ one has  $K^{-1}g(t) \leq f(t) \leq Kg(t)$. (This is the
big-Theta asymptotic notation of Bachmann--Landau.)

We fix a smooth germ $(\bC^n,0)$ and let   $(D_1,0)$ and $ (D_2,0)$ be two  curve germs in $ (\bC^n,0)$.

 If we denote by $d_o$ the outer distance $d_o(x, y) = ||x-y||_{\bC^n}$, and let $S_{\epsilon} \subset \bC^n$ be the sphere around $0$ with radius
 $\epsilon$, then there exists a rational number $q_0(D_1,D_2)$, called the {\it outer contact of $(D_1,0)$ and $(D_2,0)$},  such that
$$\big(\, \epsilon\mapsto d_o(D_1 \cap S_{\epsilon}, D_2 \cap S_{\epsilon})\,\big) = \Theta(\epsilon^{q_{o}(D_1, D_2)}).$$
Now, if $(X,0)$ is a normal surface singularity, $(X,0)\subset (\bC^n,0)$, and $(D_1,0)\cup(D_2,0)\subset (X,0)$, then
we can consider also the inner metric $d_i$ on $(X, 0)$ defined by $d_i^X(x, y) = \min_{\gamma \in \mathcal{C}} \text{length}(\gamma)$, where $\mathcal{C}$ is the space of curves
$\gamma: [0, 1] \to X$ such that $\gamma(0) = x$ and $\gamma(1) = y$.
Then  the {\it inner contact}   $q_i^{X}(D_1, D_2)\in\mathbb{Q}$  is defined  via:
$$\big(\, \epsilon\mapsto d_i^X(D_1 \cap S_{\epsilon}, D_2 \cap S_{\epsilon})\,\big) = \Theta (\epsilon^{q_i^{X}(D_1, D_2)} ).$$
Notice that $d_i^X(x, y) \geq d_o(x, y)$, which implies that $q_i^{X}(D_1, D_2) \leq q_{o}(D_1, D_2)$.

\smallskip

 Let us fix also a generic linear projection $l: (z_1, z_2, \dots, z_n) \to (z_1, z_2)$, where $(z_1, z_2, \dots, z_n)$ are local analytic
   coordinates of $(\bC^n,0)$ and $(z_1, z_2)$ are local analytic coordinates of $(\bC^2,0)$.

The inner and outer metrics coincide on the smooth germ $(\bC^2, 0)$ (that is, when
$(X,0)=(\bC^2,0)$); denote this metric by $d^{\bC^2}$ and the corresponding exponent
by $q^{\bC^2}$.
Since for any two points $x, y \in (X, 0)$ we have $d_{o}(x, y) \geq d^{\bC^2}(l(x), l(y))$,
it follows that:
$$q_i^{X}(D_1, D_2) \leq q_{o}(D_1, D_2) \leq q^{\bC^2}(l(D_1), l(D_2)).$$
Next, we recall a computaional formula for $q_0(D_1,D_2)$. For simplicity we will  assume that $(D_1,0)$ and $(D_2,0)$  have the same multiplicity $m$
(this is the case what we will need below).
 Thus, they admit
  local parametrizations of   the form  $(D_1, 0) = (g_1(t), \dots, g_n(t))$, where $g_1(t) = t^m$
and $g_i(t) = \sum_{j \geq m} a_{i, j} t^j$, and $(D_2, 0) = (h_1(t), \dots, h_n(t))$, where $h_1(t) = t^m$ and $h_i(t) = \sum_{j \geq m} b_{i, j} t^j$, $i\geq 2$.

Let $\omega$ be a primitive $m$-th root of unity, and for an integer $1 \leq l \leq m$ denote:
\begin{equation}\label{eq:Jl}
J_{\,l}  = \min  \{ j \mid \exists \ 1 \leq i \leq n, b_{i, j} \neq a_{i, j} \cdot \omega^{l j} \}.\end{equation}
Then (see e.g. \cite{neumpichoncurve}),
\begin{equation}\label{eq:q12} q_{o}(D_1, D_2) := \frac{1}{m} \cdot \max_{1 \leq l \leq m} J_l.\end{equation}

\medskip

In the next discussion we will make a special choice of the pair $(D_1,D_2)$. Let us fix $(X,0)\subset (\bC^n,0)$ and let $\tX\to X$ be a good resolution.
We fix  one of the  irreducible  exceptional divisors, $E_v$. Consider two generic smooth transversal curvettes $C_1, C_2 \subset \tX$ which meet
 $E_v$ at distinct points. Consider the projections  of these curvettes $D_1 = \pi (C_1)$ and $D_2 = \pi (C_2)$.
Clearly, they have the same multiplicity $m$, where $m$ is the  $v$-coefficient of the maximal ideal cycle.

\begin{theorem}\label{curvettesinner} \cite[Lemma 3.2]{fantini}

Consider two curvettes $C_1$ and $C_2$ of $E_v$ meeting it at generic distinct points, and let $D_1 = \pi(C_1)$ and $D_2 = \pi(C_2)$.
Then the inner rate $q_i^X(D_1 , D_2)$ between $D_1$ and $D_2$ depends only on $v$ and not on the choice of $C_1$ and $C_2$. Moreover,    $q_i^X(D_1 , D_2) \geq 1$,
$m \cdot q_i^X(D_1 , D_2)$ is an integer, and if $l: (X,0) \to (\mathbb{C}^2,0)$ is a generic projection with respect to $\pi$, then 
\[ q_i^X(D_1 , D_2) = q^{\mathbb{C}^2} \left(l(D_1), l(D_2) \right). \]
In particular, $q_i^{X}(D_1, D_2) = q_{o}(D_1, D_2) = q^{\bC^2}(l(D_1), l(D_2))$.

Furthermore, the rational number $q_i^{X}(D_1, D_2) $
is independent of the local analytic embedding of $(X,0)$ into some $(\bC^n,0)$,
in particular it is also independent of the choice of the local analytic coordinates
of $(\bC^n,0)$ (hence of the choice of the corresponding associated Hermitian metric too).
\end{theorem}

\begin{definition}
In the sequel the common value
 $q_i^{X}(D_1, D_2) = q_{o}(D_1, D_2) = q^{\bC^2}(l(D_1), l(D_2))$ associated with $E_v$ will be called  the {\it   inner rate of  $v\in\mathcal{V}$}  and will be denoted by $q_v$.
 \end{definition}

%
\subsection{Simultaneous parametrization}\label{ss:simult}
Above, in the formula (\ref{eq:q12}),
 the parametrizations $t\mapsto g(t)$ and $t\mapsto h(t)$ of $(D_1,0)$ and $(D_2,0)$ respectively have no  other compatibility relations  than
$g_1(t)=h_1(t)=t^m$. Any of them can be reparametrized by $t\mapsto \omega \cdot t$, and this reparametrization
ambiguity is resolved by considering all the integers $\{J_l\}_{l=1}^{m}$
and the  expression $q_0=\frac{1}{m}\max_l\{J_l\}$.
We prove that if we make the parametrizations more compatibly, then $\max_l\{J_l\}$ will be realized by $J_m$.
(That is, in the comparison of the coefficients  $a_{i,j}$ and $b_{i,j}$ no root of unity will be needed anymore.)
 This  construction will be used in the proof of Theorem \ref{ciclicinner}.

Fix  a generic point $p$  of $E_v$, let $U_p$ be a local neighborhood of $p$ in $\tX$ . Let  $L$ be a fixed  generic linear function, and choose local coordinates
$(\alpha,\beta)$  in $U_p$    such that $\pi^*L=\alpha^m$ in $U_p$. Let us define the cut $C_1$ in $U_p$ given by $\beta=0$ and $C_2$ by
$\beta=\beta_0$,  a nonzero constant with sufficiently small absolute value (which will be determined below). Set $(D_1,0)=\pi(C_1)$, $(D_2,0)=\pi(C_2)$. Let
$z_1=L, z_2, \ldots, z_n$ be  generic coordinates in $(\bC^n,0)$, their pullbacks restricted to $U_p$ has the form
$z_1=L=\alpha ^m, \ \  z_i=\sum_{j\geq m}a_{i,j}(\beta)\alpha^j$ (for $i\geq 2$), where $a_{i,j}(\beta)\in \bC\{\beta\}$.
Then the parametrization of $\pi (\{ \alpha=t, \, \beta=\beta_0 \})$ has the form
\begin{equation}
\label{eq:para}z_1=t^m, \ \  z_i=\sum_{j\geq m}a_{i,j}(\beta)|_{\beta=\beta_0}t^j
\ \  \  (i\geq 2).
\end{equation}

Let $J_m':= \min  \{ j \mid \exists \ 2 \leq i \leq n, a_{i, j}(\beta) \not\equiv a_{i, j}(0). \}$
Let $(i_0,j_0)$ be a pair which realizes the minimum in $J_m'$. Then there exists $\epsilon >0$ sufficiently small such that
$a_{i_0,j_0}(\beta)\not= a_{i_0, j_0}(0)\cdot \omega ^{lj_0}$ for any $l$ and $0<|\beta|\leq \epsilon$. Therefore, for $0<|\beta_0|\leq \epsilon$,
and for the two curves $(D_1,0)$ ad $(D_2,0)$ defined above,  the $J_l$ values satisfy: $J_m=J_m'$, and $J_l\leq J_m'$, hence
$q_0(D_1,D_2)=\frac{1}{m}\cdot J_m$.

\section{The first   $q_v$ expressions}

In this section $(X,0)$ is an arbitrary normal surface singularity with rational homology sphere link, as in \ref{ss:notation}.

\subsection{Preliminary remarks regarding base points}

Let us consider the situation of  Theorem \ref{curvettesinner} with $C_1,\, C_2\in \tX$, two
generic curvettes of $E_v$. In order to compute their inner rate
we wish to separate them by a section of a certain line bundle $\mathcal{L}\in {\rm Pic}(\tX)$. This means that we would like to construct a section $s$ of $\mathcal{L}$  such that the zero set of $s$ intersects $C_1$ but not $C_2$ (along $E$). In order to realize this (since $C_1$ and $C_2$ are generic), the intersection point of the zero set of $s$ with $E_v$ should be movable along $E_v$, that is, $\mathcal{L}$ must be base point free along $E_v$.
For any  (a priori identified)  line bundle, we can eliminate its base points by finitely many (iterated)  blow ups (where the centers are determined by the base points).
On the other hand, if we do not know a priori which line bundle will appear along our  construction, then we have to develop different `elimination techniques'.

\medskip
In this subsection we will exemplify both cases.

\medskip

In the first case our line bundle  is (the easily identified)  $\mathcal{L}=\calO_{\tX}(-Z_{max})$.
Then,  it is known  (see e.g. \cite{fantini}),  that if $(Z_{max}, E_v)<0$ and
$\mathcal{L}$ has no base point along $E_v$ then $q_v=1$.
Indeed,  let $l$ be the generic projection  $(X,0)\subset (\bC^n,0)\to (\bC^2,0)$
 and set
$Z_{max,v}=m$.
Then,   under the above  assumptions, the family of curvettes along $E_v$ are projected via $l$ in a family is plane curve singularities with multiplicity $m$ and different tangent cones. Then, via the definition from subsection \ref{ss:1.7}, we get  that $q_v= q^{\bC^2}(l(D_1), l(D_2))=1$.

In fact, whenever  $\tX$ factors through the blow up of the maximal ideal, the  line bundle $\calO_{\tX}(-Z_{max})$ is base point free.
Therefore,
 for any such $\tX$,
\begin{equation}\label{eq:q1}
\mbox{  for any $E_v$ which satisfies  $(Z_{max}, E_v)<0$ one has  $q_v=1$,}\end{equation}
 independently of the value $Z_{max, v}$.

 Hence, we only have to handle the case $(Z_{max}, E_v)=0$ (see e.g. Proposition \ref{innerrateone}).

In the second case, instead of proving the base point freeness of a line bundle (or instead of making some  blow ups), we will change the line bundle to another one,  which is base point free along $E_v$.

\begin{lemma}\label{nobasepoint}
Let $v \in \calv$ be an arbitrary vertex such that $(Z_{max}, E_v) = 0$  and  $Z_{max, v} = 1$. Assume that  there exists a cycle $l\in L$ with   $(l, E_v) < 0$
and which can be represented as  ${\rm div}_E(\pi^*f) = l $ for a certain function
 $f \in \calO_{(X,0)}$.  
 Then there exists another
function $g\in \calO_{(X,0)}$ such that the cycle $l':={\rm div}_E(\pi^*g)$ satisfies
 $l'_v = l_v$ and $(l',E_v)<0$, and the line bundle $\calO_{\tX}(-l')$ has no base points along $E_{v}^\circ:= E_v\setminus \cup_{u\not=v}E_u$.
\end{lemma}
\begin{proof}
Define the cycle $l'= \min\{l_v \cdot Z_{max}, l\}$.
Consider a generic linear function $L \in H^0(\calO_{\tX}(-Z_{max}))$ with $\text{div}_E(L)=Z_{max}$.
Then for any point $p \in E_{v}^\circ$ and   for a generic value $c \in \bC$, the function  $g_c= c \cdot L^{l_v} + f$ has no strict transform vanishing at $p$.
\end{proof}

\subsection{The  $(Z_{max}, E_v) = 0$  and  $Z_{max,v}=1$ case}


\begin{theorem}\label{innerrateone} Let $(X,0)$ be a normal surface singularity with rational homology sphere link  and let $\pi:\tX\to X$ be  one of its resolutions.
Let $E_v$ be an irreducible exceptional divisor
such that  $(Z_{max}, E_v) = 0$  and  $Z_{max, v} = 1$.  Then $q_v={\rm Min}_v$, where
\begin{equation}\label{eq:minv}{\rm Min}_v := \min \{\, k \,\mid\, \exists f \in \calO_{(X,0)}, \ {\rm div}_E(\pi^*f)_v = k, \
{\rm St}_{\tilde{X}}(f) \cap E_{v}^\circ \neq  \emptyset \}.\end{equation}
\end{theorem}
\begin{remark}\label{rem:ab}
 (a) If $\pi$ factors through the blow up of the maximal ideal then,
in fact, $q_v={\rm Min}_v$ is true under the assumption  $Z_{max, v} = 1$ even if   $(Z_{max}, E_v) <0$. Indeed,  in that case $q_v=1$ by (\ref{eq:q1}) and ${\rm Min}_v=1$
since in  (\ref{eq:minv}) the minimum is realized by  $f=L$.

(b) If  $(Z_{max}, E_v) = 0$  and  $Z_{max, v} = 1$ (as in Theorem \ref{innerrateone}), then $q_v\geq 2$. Indeed,   the existence of a function
$ f \in \calO_{(X,0)} $ with  $ {\rm div}_E(\pi^*f)_v = 1$ and
${\rm St}_{\tilde{X}}(f) \cap E_{v}^\circ \neq  \emptyset$ would imply  $(Z_{max}, E_v) < 0$.
\end{remark}

\begin{proof}[Proof of  Theorem \ref{innerrateone}]
Consider two generic curvettes $C_1,   C_2$ of the exceptional divisor $E_v$, and set  $(D_1,0) = \pi (C_1)$ and $(D_2,0) = \pi (C_2)$ as above. From the assumptions follows that the pullback of the generic linear section, as a divisor in $\tX$, has intersection multiplicity one with $C_a \  (a \in \{1, 2\})$, hence the intersection multiplicity at 0 of the generic linear section of
$(\bC^n,0)$ with $(D_a,0) \ (a \in \{1, 2\})$ is one as well. Hence, both $(D_1,0)$ and $(D_2,0)$ are smooth germs.

\smallskip

We claim that there exists a smooth $2$-dimensional germ  $(F,0)\subset (\bC^n,0)$  such that $(D_1,0)\cup(D_2,0) \subset (F,0)$. Indeed, assume first that the tangent spaces $T_0(D_1) \neq T_0(D_2)$. Then, if $(\bC,0) \ni t_a \mapsto f^{(a)}(t_a)\in (\bC^n,0)$ parameterizes $(D_a,0)$ then $(t_1,t_2)\mapsto f^{(1)}(t_1)+f^{(2)}(t_2)$ parameterizes a smooth embedded surface germ containing both curve germs.

Otherwise, assume that  $T_0(D_1) = T_0(D_2)$.
Consider coordinates $x_1, x_2, \dots, x_n$ on $\bC^n$ such that $(D_1, 0)$ is parameterized locally
by $(t, 0, \dots, 0)$ and  $(D_2,0)$  by $ \text{im}(f) \subset \bC^n$, where $f_1(s) = s$ and $f_i(s) = \sum_{j \geq 2} b_{i, j} s^j$ ($2\leq i\leq n$)
for some complex numbers $b_{i, j} \in \bC$.
 Let $k$ be the smallest index such that there exists an index $i \in \{2, \dots, n\}$ with $b_{i, k} \neq 0$.
Then we define a local embedding
 $F: (\bC^2, 0) \to (\bC^n, 0)$
by components
$F(t, s)_1 = t$ and $F(t, s)_i = s \cdot \sum_{j \geq k} b_{i, j} t^{j-k}$ for $i\geq 2$.
Clearly $(D_1,0)$ is parameterized by $F(t, 0)$   and $(D_2,0)$ by  $F(t, t^k)$, and $F$ is a local embedding.

\medskip

Next, consider a generic projection $l: (\bC^n,0) \to (\bC^2,0)$. Then
 $l$ induces an isomorphism between the triples of germs $(F, D_1, D_2)$ and $(\bC^2, l(D_1), l(D_2))$.

For the pair of smooth curves $(l(D_1), l(D_2))$, it is known from \cite{neumpichoncurve} that $q^{\bC^2}(l(D_1), l(D_2))$ equals  the minimal number $q$ of  consecutive  blow-ups of
$(\bC^2,0)$ needed to make the strict transforms of $l(D_1)$ and $l(D_2)$ disjoint.
  This fact reformulated states that   $q$ equals also with  the minimal number of blow ups (creating exceptional curves
  $E_1,\ldots,  E_q$, where $E_q$ is the last one) such that there exists a function
  $g \in (\calO_{\bC^2},0)$ having a strict transform
    $D$ along $E_{q}$ which intersects  the strict transform of   $l(D_1)$
     but not the strict transform of $ l(D_2)$.
     In other words, $q_v={\rm Min}_{\bC^2}$, where
     $${\rm Min}_{\bC^2}:= \min \{ \, k \, |\, \exists g\in \calO_{(\bC^2,0)}, \
   k = \min\{\text{ord}_{l(D_1)}(g), \text{ord}_{l(D_2)}(g)\},\
    \text{ord}_{l(D_1)}(g)\not= \text{ord}_{l(D_2)}(g) \}.$$
    Here $\text{ord}_{l(D_a)}(g)$ is the intersection multiplicity at 0 of $l(D_a)$ with $g$,
    that is, it equals  $\text{ord}_t(g(h(t)))$, where $h(t)$ parameterizes $l(D_a)$, $a\in\{ 1,2\}$.

\medskip

Since the projection $l|_F: (F,0) \to (\bC^2,0)$ is an isomorphism, we also have
$q_v={\rm Min}_{F}$, where
 \begin{equation}\label{eq:1}
 {\rm Min}_{F}:= \min \{ \, k \, |\, \exists g\in \calO_{(F,0)}, \
   k = \min\{\text{ord}_{D_1}(g), \text{ord}_{D_2}(g)\},\
    \text{ord}_{D_1}(g)\not= \text{ord}_{D_2}(g) \}.\end{equation}


Since $\calO_{(\bC^n,0)}\to\calO_{(F,0)}$ is surjective, the germ $g$ from (\ref{eq:1}) can be lifted to $\calO_{(\bC^n,0)}$. Moreover, since $(D_1,0)\cup (D_2,0)\subset (X,0)$, the right hand side of (\ref{eq:1}) equals with a similar expression defined for functions from $\calO_{(X,0)}$. That is, $q_v={\rm Min}_X$, where
\begin{equation}\label{eq:2}
 {\rm Min}_{X}:= \min \{ \, k \, |\, \exists g\in \calO_{(X,0)}, \
   k = \min\{\text{ord}_{D_1}(g), \text{ord}_{D_2}(g)\},\
    \text{ord}_{D_1}(g)\not= \text{ord}_{D_2}(g) \}.\end{equation}

Finally, we  prove that ${\rm Min}_X={\rm Min}_v$ (where ${\rm Min}_v$ is the expression from
the statement of the theorem).

Consider  some $g$ which realizes ${\rm Min}_X$ in (\ref{eq:2}), and consider the resolution $\pi:\tX\to X$ as above. Since  $\text{ord}_{D_1}(g)\not= \text{ord}_{D_2}(g)$, the strict transform of the divisor of $g$ via $\pi$  must intersect at least one of the curvettes $C_a$, hence
${\rm St}_{\tilde{X}}(g) \cap E_{v}^\circ  \neq  \emptyset $ and also $({\rm div}_E(\pi^*g))_v\leq {\rm Min}_X$
(since ${\rm Min}_X$ is the sum of $({\rm div}_E(\pi^*g))_v$ with the intersection multiplicity of $C_a$ which realizes the minimum with the strict transform of $g$).
But then, by the definition of ${\rm Min}_v$, we also have ${\rm Min}_v\leq ({\rm div}_E(\pi^*g))_v$. Hence ${\rm Min}_v\leq {\rm Min}_X$.

Conversely, consider some $f$ which realizes ${\rm Min}_v$, i.e. ${\rm div}(\pi^*f)_v={\rm Min}_v$ and ${\rm St}_{\tilde{X}}(f) \cap E_{v}^\circ
 \neq  \emptyset$. Then by Lemma \ref{nobasepoint}
 we can assume additionally that
$\calO_{\tX}(-{\rm div}_E(\pi^*f))$ has no base points on $E_{v}^\circ$. Then for generic $C_1$ and $C_2$ we can find in this linear system another function $f' $ with the following properties:
 ${\rm div}(\pi^*f')_v={\rm Min}_v$,  the strict transform of $f'$ intersects $C_1$ but not $C_2$, in particular, ${\rm Min}_v={\rm ord}_{D_2}(f')< {\rm ord}_{D_1}(f')$.
 So $f'$  can be considered as one of the germs from the right hand side of (\ref{eq:2}) (whose $k$ is ${\rm Min}_v$). Hence ${\rm Min}_v\geq {\rm Min}_X$.

 All these together show $q_v={\rm Min}_X={\rm Min}_v$.
 \end{proof}

\subsection{$q_v$  and  the analytic semigroups associated with the resolutions}\label{ss:AS}

Recall that for a fixed resolution $\pi$, the analytic semigroup associated with $\pi$ is defined as
$\calS_{an}(\pi):=\{ {\rm div}_E(\pi^*f)\, :\, f\in \calO_{(X,0)}\}\subset L(\pi)$. If the resolution $\pi'$ dominates $\pi$, that is, if there exists a modification $\phi:\tX'\to \tX$
such that $\pi'=\pi\circ \phi$, then the natural homological projection $L(\pi')=H_2(\tX',\Z)\to H_2(\tX,\Z)=L(\pi)$ sends $\calS_{an}(\pi')$ onto
$\calS_{an}(\pi)$.
In particular, $\calS_{an}(\pi')$ contains more information than  $\calS_{an}(\pi)$ (e.g. about those strict transforms of $\pi$ which  pass the singular points of $E$).

Note also that if $\mathcal{T}$ is the dual graph of $\pi$ and $v$ is one of its vertices, then $\calS_{an}(\pi,v):=\{s\in \calS_{an}(\pi)\, :\,
(s, E_v)<0\}$ is a sub-semigroup of $\calS_{an}(\pi)$ (with $0\not\in \mathcal{S}_{an}(\pi,v)$). Thus, its projection to the $E_v$-coordinate
\begin{equation}\label{eq:san}
\overline{\calS}_{an}(\pi,v):= \{ k \, : \,  \exists s\in \mathcal{S}_{an}(\pi), \ s_v=k, \ (s, E_v)<0\}\end{equation}
is a numerical sub-semigroup in $\mathbb{N}$. Observe that if $s\in \mathcal{S}_{an}\setminus \{0\}$ then all the $E_v$-coordinates of $s$ are strict positive.
Hence, $0\not \in \overline{\calS}_{an}(\pi,v)$.

We wish to connect $\overline{\calS}_{an}(\pi,v)$  vith $q_v$ via the expression ${\rm Min}_v$ from (\ref{eq:minv}).

First  observe that for $s={\rm div}_E(\pi^*f )$, if ${\rm St}_{\tX}(f)\cap E^\circ _v\not=\emptyset$  then $(s, E_v)<0$ holds, however, the converse is not
necessarily true. Therefore, in order to connect (\ref{eq:san}) with (\ref{eq:minv}) we need to perform an additional construction. We start with the following lemmas.

\begin{lemma}\label{lem:free} (See \cite[Theorem 6.5.7]{NBook} Fix a resolution $\pi$. For any vertex $u$ of $\mathcal{T}$
 let $\delta_u$ denote the number of adjacent vertices of $u$ in $\mathcal{T}$.  Let $D\in L$  be a cycle such that  $(D,E_u)\leq -\delta_u$  for any $u$. Then the line bundle
$\calO_{\tX}(-D)$ is base point free and the zero set of its generic section intersects (transversely)
 each irreducible exceptional curve $E_u$.
\end{lemma}
This implies that for any fixed vertex $v$, there exists $g\in \calO_{(X,0)}$ such that ${\rm St}_{\tX}(g)\cap E^\circ _v\not=\emptyset$.
(In particular, ${\rm Min}_v$ in (\ref{eq:minv}) is well-defined.) Fix such a function $g$ and set $N_v:= {\rm div}_E(\pi^*g)_v$. In particular, in the determination  of 
${\rm Min}_v$ in (\ref{eq:minv}) we have to test only those functions $f$  which satisfy  $ {\rm div}_E(\pi^*f)_v\leq N_v$. 
Next, fix some $r_v\in \mathbb{Q}_{>0}$ which satisfies Lemma \ref{lem:iv} for $\pi$ and $v$. Furthermore, take any integer $b_v> r_v\cdot N_v$.

Then we modify the resolution $\pi$ by blowing up each $\{E_v\cap E_u\}_{u\not= v}$ $b_v$ times. By this we mean that  we blow up $E_v\cap E_u$,
then we blow up the intersection of the strict transform of $E_v$ with the newly created curve, and we continue in this way. 
Let $\pi'$ the newly created resolution, $E'=\pi'^{-1}(0)$  and $E_v'$ the strict transform of $E_v$.  

\begin{lemma}\label{lem:divst}
If for some function $f\in\mathcal{O}_{(X,0)}$ one has 
${\rm div}_{E'}(\pi'^*f)_v\leq N_v$ and 
$(\, {\rm div}_{E'}(\pi'^*f), E_v'\, )_{\mathcal{T}'}<0$ then ${\rm St}_{\tX'}(f)\cap E_v'^\circ\not=\emptyset$.
\end{lemma}
\begin{proof} Since $(\, {\rm div}_{E'}(\pi'^*f), E_v'\, )_{\mathcal{T}'}<0$ the strict transform ${\rm St}_{\tX'}(f)$ must intersect $E_v'$ somewhere.
Assume that  ${\rm St}_{\tX'}(f)\cap E_v'^\circ=\emptyset$. Then ${\rm St}_{\tX'}(f)$  intersects $E_v'$ in an intersection point with some of its neighbours.
But then,  $E_v$ in $\tX$ must have a neighbour $E_u$ such that  ${\rm St}_{\tX}(f)$ passes at $E_v\cap E_u$ and 
in fact  $({\rm St}_{\tX}(f), E_v)\geq b_v$. But then (cf. Lemma \ref{lem:iv})
 $N_v\geq  {\rm div}_{E'}(\pi'^*f)_v={\rm div}_{E}(\pi^*f)_v\geq ({\rm St}_{\tX}(f), E_v)/r_v \geq b_v/r_v> N_v$, a contradiction. 
\end{proof}

\begin{corollary}\label{cor:case1} Assume that for some resolution $\pi$ we have $(Z_{max}, E_v)=0$ and $Z_{max,v}=1$. 
Let $\pi'$ be the resolution constructed above from $\pi$. Then $q_v= \min \, \overline{\calS}_{an}(\pi',v)$.
In particular, $q_v$ can be computed from the analytic semigroup of one of its resolutions which can be completely identified 
(that is, if we blow up the intersection points 
$\{E_v\cap E_u\}_u$ in $\tX$ sufficiently many times, where the number of blow ups can be determined). 
\end{corollary}
\begin{proof} Combine Theorem \ref{innerrateone} with Lemma \ref{lem:divst}. 
\end{proof}
Recall that if $(X,0)$ is rational, then any resolution is good, and for any resolution $Z_{max}=Z_{min}$ (where 
 $Z_{min}=\min \{ \mathcal{S}\setminus \{0\}\}$ is the Artin minimal/fundamental cycle), 
  $\calO_{\tX}(-Z_{max})$ is base point free and 
$\calS_{an}=\calS$, see \cite[7.1.A]{NBook}. 
\begin{corollary}\label{cor:rat}
Assume that $(X,0)$ is rational and let $\pi$ any  of its  resolutions. Fix a vertex $v$.  Then $q_v$ is topological  whenever one of the following conditions hold: 
  (i) $(Z_{min}, E_v)<0$, or (ii) $(Z_{min}, E_v)=0$ and $Z_{min,v}=1$.  (We will add  a third case later in subsection \ref{ss:5.1}.)
\end{corollary}

\subsection{Examples} \label{ss:examples}\
In the following examples we show that in the non-rational case Corollary \ref{cor:rat} is not necessarily true (even if the link is integral  homology sphere, or even for hypersurface Newton non-degenerate singularities, or Milnor number constant deformations). 

\begin{example}\label{ex:11} Consider the following graph.

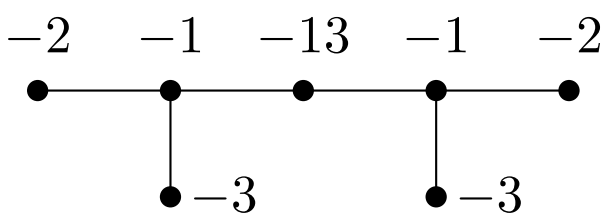
\begin{picture}(300,45)(100,0)
\put(125,25){\circle*{4}} \put(150,25){\circle*{4}}
\put(175,25){\circle*{4}} \put(200,25){\circle*{4}}
\put(225,25){\circle*{4}} \put(150,5){\circle*{4}}
\put(200,5){\circle*{4}} \put(125,25){\line(1,0){100}}
\put(150,25){\line(0,-1){20}} \put(200,25){\line(0,-1){20}}
\put(125,35){\makebox(0,0){$-2$}}
\put(150,35){\makebox(0,0){$-1$}}
\put(175,35){\makebox(0,0){$-13$}}
\put(200,35){\makebox(0,0){$-1$}}
\put(225,35){\makebox(0,0){$-2$}} \put(160,5){\makebox(0,0){$-3$}}
\put(210,5){\makebox(0,0){$-3$}}
\end{picture}

 The corresponding link is an integral homology sphere.
This graph supports several  analytic types, for their classification see \cite{okuma}. The possible values of the geometric genus is $2\leq p_g\leq 3$.
E.g., for the generic analytic structure $p_g=2$. Here we will treat the two cases when   $p_g=3$.  For more details see \cite{okuma}. 
Let $v$ be the vertex with decoration $(-13)$.

One of the analytic types is a (splice type) complete intersection, with equations $z_1^2z_2+z_3^2+z_4^3=z_1^3+z_2^2+z_4^2z_3=0$ in $(\mathbb{C}^4,0)$.
Its maximal ideal cycle $Z_{max}$ equals $2 Z_{min}$.
The strict transform of $L$ has two components, both intersect  $E_v^\circ$.  Since $Z_{min,v}=1$, $(Z_{max}, {\rm St}_{\tX}(L))=4$, which is 
also the multiplicity of $(X,0)$, cf. \cite{okuma},  $\calO_{\tX}(-Z_{max})$ is base point free. Hence, by (\ref{eq:q1}), $q_v=1$.

Next, we consider a Kodaira singularity type supported on $\mathcal{T}$. It is not Gorenstein. 
 In this case $Z_{max}=Z_{min}$. The strict transform ${\rm St}_{\tX}(L)$ 
of $L$ is irreducible, and intersects $E_v$ transversely along $E_v^\circ$.  But this intersection point is a base point of $\mathcal{O}_{\tX}(-Z_{max})$
(since the multiplicity equals 3,  while  $-Z_{max}^2=-Z_{min}^2=1$). In particular, in order to get the resolution  which factors through the blow up of the maximal ideal, we have to blow up (twice) the intersection point ${\rm St}_{\tX}(L)\cap E_v$. In this resolution $Z_{max, v}=1$ and $(Z_{max}, E_v)=0$. Hence 
$q_v>1$ by Remark \ref{rem:ab}. 
\end{example}

\begin{example}\label{ex:22}
We consider the following one-parameter family of hypersurface singularities in $(\mathbb{C}^3,0)$ given by 
$f_t(x,y,z)=x^5+y^3+z^{10}+x^4z+tx^3y=0$, where $t\in (\mathbb{C},0)$ is the deformation parameter. 
This is the analogue of the classical Brian\c{c}on-Speder type deformation \cite{BS}, 
however in our present case we have additionally that the stable link is a 
rational homology sphere. Along the deformation   the link, the embedded topological type 
(hence the Milnor fibration, the geometric genus) are stable, it admits a weak simultaneous resolution, but does not admit a strong simultaneous resolution, or equivalently, it 
  is not 
Whitney equisingular.   The maximal ideal cycle jumps (see discussion below), and the Milnor fiber of the generic linear section also jumps. 
In the next discussion  we show  that the Lipschitz inner rate of a certain vertex jumps as well. 

For both cases $t=0$ and $t\not=0$, the hypersurface singularity $f_t$ has non-degenerate Newton principal part. In particular, the embedded 
minimal good resolution graph 
of $(X_t,0)=(f_t=0,0)$, together with any of the coordinate functions $x,y,z $, can be computed by Oka's Algorithm \cite{oka}, see also \cite[5.5.B]{NBook}
or \cite{NSi}. The very same vertex $v$ is marked in both graphs.  

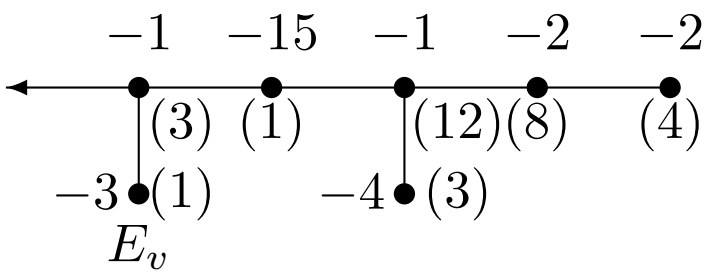
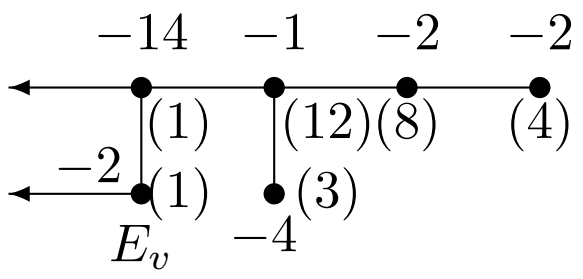
\vspace{2mm}

\begin{picture}(300,55)(80,-10)
\put(150,25){\circle*{4}} \put(250,25){\circle*{4}}
\put(175,25){\circle*{4}} \put(200,25){\circle*{4}}
\put(225,25){\circle*{4}} \put(150,5){\circle*{4}}
\put(200,5){\circle*{4}} \put(150,25){\line(1,0){100}}
\put(150,25){\line(0,-1){20}} \put(200,25){\line(0,-1){20}}
\put(150,25){\vector(-1,0){25}}
\put(250,35){\makebox(0,0){$-2$}}
\put(150,35){\makebox(0,0){$-1$}}
\put(175,35){\makebox(0,0){$-15$}}
\put(200,35){\makebox(0,0){$-1$}}
\put(225,35){\makebox(0,0){$-2$}}
 \put(140,5){\makebox(0,0){$-3$}}
\put(190,5){\makebox(0,0){$-4$}}

\put(250,18){\makebox(0,0){$(4)$}}
\put(158,18){\makebox(0,0){$(3)$}}
\put(175,18){\makebox(0,0){$(1)$}}
\put(210,18){\makebox(0,0){$(12)$}}
\put(225,18){\makebox(0,0){$(8)$}}
\put(158,5){\makebox(0,0){$(1)$}}
\put(150,-5){\makebox(0,0){$E_v$}}
\put(210,5){\makebox(0,0){$(3)$}}

\put(105,15){\makebox(0,0){$(t=0):$}}\put(315,15){\makebox(0,0){$(t\not=0):$}}

 \put(450,25){\circle*{4}}
\put(375,25){\circle*{4}} \put(400,25){\circle*{4}}
\put(425,25){\circle*{4}} \put(375,5){\circle*{4}}
\put(400,5){\circle*{4}} \put(375,25){\line(1,0){75}}
\put(375,25){\line(0,-1){20}} \put(400,25){\line(0,-1){20}}
\put(375,25){\vector(-1,0){25}}\put(375,5){\vector(-1,0){25}}
\put(450,35){\makebox(0,0){$-2$}}
\put(375,35){\makebox(0,0){$-14$}}
\put(400,35){\makebox(0,0){$-1$}}
\put(425,35){\makebox(0,0){$-2$}}
 \put(365,10){\makebox(0,0){$-2$}}
\put(398,-3){\makebox(0,0){$-4$}}

\put(450,18){\makebox(0,0){$(4)$}}
\put(382,18){\makebox(0,0){$(1)$}}\put(382,5){\makebox(0,0){$(1)$}}
\put(410,18){\makebox(0,0){$(12)$}}
\put(425,18){\makebox(0,0){$(8)$}}
\put(375,-5){\makebox(0,0){$E_v$}}
\put(410,5){\makebox(0,0){$(3)$}}

\end{picture}

\vspace{2mm}

In the $t=0$ case   $Z_{max} ={\rm div}_E(z)$ (see the  multiplicities in $(\ )$ on   left diagram, the arrow corresponds to ${\rm St}_{\tX}(z)$).
In particular, $Z_{max, v}=1$ and  $(Z_{max}, E_v)=0$. On the other hand, the function $y$ has multiplicity 2 along $E_v$ and ${\rm St}_{\tX}(y)\cap E_v^\circ\not=\emptyset$.
Hence $q_v=2$  by Theorem \ref{innerrateone}. 

If $t\not=0$ then again $Z_{max} ={\rm div}_E(z)$, $Z_{max, v}=1$ but now  $(Z_{max}, E_v)<0$. On the other hand, the function $x$ also has multiplicity 1 along $E_v$
and 
 ${\rm St}_{\tX}(x)\cap E_v=\emptyset$. Hence $Z_{max}$ has no base point along $E_v$. Therefore, $q_v=1$. 

\end{example}

\section{More $q_v$ expressions. The covering case}\label{s:COV}

\subsection{The general $(Z_{max}, E_v)=0$ and  $Z_{max, v} > 1$ case} \label{ss:cover} \

Consider a generic linear function $L\in \calO_{(\mathbb{C}^n, 0)}$,  its restriction
to  $ \calO_{(X, 0)}$,  still denoted by $L$, and, as above,  $Z_{max}:=\sum_{u \in \mathcal{V}} l_u E_u := {\rm div}_E(\pi^*L)$.
We assume that the  resolution  $\tX$ factors through the blow-up of the maximal ideal,
 hence the line bundle $\calO_{\tX}(-Z_{max})$ has no base points.
Denote  the irreducible components of $L^{-1}(0)\subset X$ by $\cup _{1 \leq i \leq M} W_i$, where $M =- (Z_{max}, E)$.

We fix a vertex $v$ such that  $(Z_{max}, E_v)=0$ and  $m:=l_v > 1$.   In the following
we wish to reduce this  case to the `$Z_{max, v} = 1$ case'  via a cyclic branched covering.

\medskip

Let $X^{(m)}$ be the fiber product of $L:(X,0)\to (\bC,0)$ and $(\bC,0)\to (\bC,0)$,  $z\mapsto z^m$. That is, $X^{(m)}:= \{(x,z)\in (X\times \bC,0\times 0)\,:\, L(x)=z^m\}$. Denote the embedding $X^{(m)}\hookrightarrow X\times \bC $ by $\iota$ and the projection $ X\times \bC  \to X$ by $pr$. Then the  map $f_m = pr \circ \iota: X^{(m)} \to X$ is a cyclic covering
ramified along $L^{-1}(0)$ of degree $m$. We know that $X^{(m)}$ has an isolated singularity because $L^{-1}(0)$ is reduced.
Furthermore, $X^{(m)}$ is Cohen-Macaulay since it is defined by a single non-zero divisor  equation $L(x) = z^m$ in the Cohen-Macaulay space $X\times \mathbb{C}$.
Therefore,  $X^{(m)}$ is a {\it normal} surface singularity (by Serre criterion).

Consider our fixed embedded  resolution $\pi: (\tX, \cup_{1 \leq i \leq M} V_i) \to (X, \cup_{1 \leq i \leq M} W_i)$. This means that $\pi^{-1}(\cup_{1 \leq i \leq M} W_i)=
\pi^{-1}(0)\cup \cup_{1 \leq i \leq M} V_i$ is a normal crossing divisor. Let
$\tX^{(m)}$ be the normalized  fiber product of the maps $f_m$ and $\pi$, let
$\tilde{f}_m:\tX^{(m)}\to \tX$ and $\pi_m:\tX^{(m)} \to X^{(m)}$ be the corresponding induced maps.
For combinatorial properties of $\tX^{(m)}$ and its exceptional curves see \cite{cyclic,NBook}. 

For each $1 \leq i \leq M$ set 
$V_i':=\tilde{f}_m^{-1}(V_i)$; it is  irreducible.

In general $\tX^{(m)}$ is not smooth, over singular points
of the (reduced) divisor $\pi^{-1}(0)\cup \cup_{1 \leq i \leq M} V_i \subset \tX$ we might have cyclic quotient singularities in $\tX^{(m)}$. If we resolve them then we get a resolution
$\widehat{X}^{(m)}\to X^{(m)}$ of $X^{(m)}$. 
However, above the intersection points sitting on $E_v$ we have only trivial quotient singularities. 
Indeed, they are normalizations of germs with local equation of type $\alpha^m\beta^n=z^m$, which have ${\rm gcd}(m,n)$ irreducible components, all of them smooth.
Hence, in a neighborhood of  $\tilde{f}_m^{-1}(E_v)$ the space $\tX^{(m)}$ is smooth.  

Since in this note we wish to understand the rates $q_v$ associated with the exceptional divisors of $\pi$ via their generic curvettes, and (the coverings of)
 these curvettes already are present in the modification $\tX^{(m)}$, in our next discussions we will use the covering  $\tX^{(m)}\to \tX$,
or the partial resolution $\tX^{(m)}\to X^{(m)}$ instead of the resolution $\widehat{X}^{(m)}$.

The coverings $X^{(m)}\to X$ and  $ \tX^{(m)}\to \tX$ admit natural Galois
$\mathbb{Z}_m$ actions (via the multiplication $z\mapsto \omega\cdot z$ ($\omega =e^{2\pi i/m}$).
 Let  $\mathcal{T}_e = \mathcal{T} \cup \mathcal{T}_a$  be the extended resolution graph of $(\tX, \cup_{1 \leq i \leq M} V_i)$.
It has vertex set $\calv_e = \calv \cup \calv_a$,
where $\calv$ indexes the irreducible exceptional curves $\{E_u\}_u$ and $\calv_a$ the arrowhead vertices representing $\cup_{1 \leq i \leq M} V_i$. All of them carry the corresponding   multiplicity of $\pi^*L$: each $u\in\calv$ has multiplicity $Z_{max,u}$ while
each $a\in\calv_a$ has multiplicity one.
 The pullback of the function $L$ to $\tX^{(m)}$  equals
  $\pi_m^*(z^m)$. For the multiplicities of $\pi_m^*(z)$ in $\tX^{(m)}$ see \cite{cyclic} or \cite[4.1.C]{NBook}.

Similarly, let $\mathcal{T}^{(m)}_e = \mathcal{T}^{(m)} \cup \mathcal{T}^{(m)}_a$, with vertex set $\calv^{(m)}_e = \calv^{(m)} \cup \calv^{(m)}_a$,
denote the extended resolution graph of $(\tX^{(m)}, \cup_{1 \leq i \leq M} V'_i)$.
 Again,  the arrowhead vertices $\calv^{(m)}_a$ correspond to the noncompact  irreducible components  $\cup_i V'_i$. The geometric covering induces a graph covering
$\tau_m: \mathcal{T}^{(m)}_e \to \mathcal{T}_e$.

For $u\in\calv$, assume that $\tau_m^{-1}(u) = \{u_1, \dots, u_{k_u}\}$. Set also $\tilde{f}_m^{-1}(E_u)=\cup_{j=1}^{k_u} E^{(m)}_{u_j}$.
The cyclic Galois group $\mathbb{Z}_m$ acts transitively on the set $\tau_m^{-1}(u)$. In particular, we have the equality of the inner rates associated with
 $ \{E^{(m)}_{u_j}\}_j$:
 $q^{(m)}_{u_1}=\cdots=q^{(m)}_{u_{k_u}}$.

\bekezdes\label{bek:mult_szamolas}
We recall from \cite{cyclic,NBook} the following numerical properties of the ramified cyclic covering $\tilde{f}_m$
regarding  the behavior of the vanishing orders of the functions along the covering.

Let $p$ be a generic point on $E_u$, $U_p$ a small neighbourhood of $p$ in
$\tX$  with local coordinated $(\alpha, \beta)$ such that $U_p\cap E_u=\{\alpha=0\}$. Let $l_u=
Z_{max,u}$ be the vanishing order of $\pi^*L$ along $E_u$ and $\pi^*L=\alpha^{l_u}$ in $U_p$.  Then in the fiber product of $\pi$ and $f_m$  above $U_p$  we have the space $\{(\alpha,\beta,z)\,:\, \alpha^{l_u}=z^m\} $. This has $d_u={\rm gcd}(l_u,m)$ irreducible components of type
$\{(\alpha,\beta,z)\,:\, \alpha^{l_u/d_u}=\xi\cdot z^{m/d_u}\} $ $(\xi^{d_u}=1)$.
Their normalizations are smooth with local coordinates $(w,\beta)$,  where $\alpha=w^{m/d_u}$
and $z=c\cdot w^{l_u/d_u}$ (for some $c\in\bC^*$). Hence the vanishing order of $z$ along any $E^{(m)}_{u_j}$ is $l_u/d_u=l_u/(l_u,m)$. This, in the case of $v=u$, it equals one.

Consider any $g\in\mathcal{O}_{(X,0)}$  such that $\pi^*g$ has vanishing order $n_u$ along $E_u$.
Then the local identity  $\pi^*(g)|_{U_p}=\alpha^{n_u}$ lifts to $\alpha^{n_u}=w^{n_u\cdot m/d_u}$, hence the vanishing order of the pullback of $g$ above $U_p$ in $\tX^{(m)}$ is
$n_u\cdot m/(l_u,m)$.

Since ${\rm div}_E(\pi^*g)\geq {\rm div}_E(\pi^*L)$ for any $g\in \mathcal{O}_{(X,0)}$
(i.e. $n_u\geq l_u$),
 we also have $n_um/(l_u,m)\geq l_um/(l_u,m)$, or ${\rm div}_{E^{(m)}}(\pi_m^*f_m^{*}g) \geq
 {\rm div}_{E^{(m)}}(\pi_m^*f_m^{*}L)$, hence
${\rm div}_{E^{(m)}}(\pi_m^*f_m^{*}g)\geq
m \cdot {\rm div}_{E^{(m)}}(\pi_m^{*}z)>{\rm div}_{E^{(m)}}(\pi_m^{*}z)$.  Note also
 that the functions $\{\{f_m^*g\}_{g\in \mathcal{O}_{(X,0)}}, z \}$ generate $\mathcal{O}_{(X^{(m)},0)}$ as a ring. Thus:
\begin{equation}\label{zmaxm}
{\rm div}_{E^{(m)}}(z) \ \mbox{ is the maximal ideal cycle $Z_{max}^{(m)}$ \ in $\tX^{(m)}$.}
\end{equation}
We prove next the main comparision theorem connecting the inner rates of the covering and the original singularity:

\begin{theorem}\label{ciclicinner}
The inner rates \ $q^{(m)}_{v_1}$ and $q_v$ satisfy \
$q^{(m)}_{v_1} = m \cdot q_v$.
\end{theorem}
\begin{proof}
We can assume that the
 generic coordinate system on $(\bC^n, 0)$  and $L$ was chosen in such a way that  $x_1 = L$.
We will parametrize $(D_1,0)$ and $(D_2,0)$ `simultaneously'  as in subsection \ref{ss:simult}, whose notation we will follow.

Let us analyze next the covering $\cup_j E^{(m)}_{v_j}\to E_v$. If $p$ is a generic point of $E_v$ and $U_p$ one of its small neighbourhoods with local coordinates
 $(\alpha,\beta)$ as above, then the covering over $U_p$ is given by $\alpha^m=z^m$, hence $p$ has $m$ preimages, and
$U_p$ has also $m$ preimages $U^{(m)}_{p,1}, \dots, U^{(m)}_{p,m}$, indexed in such a way that
the Galois action    $z\mapsto z\cdot e^{2\pi i/m}$ maps $U^{(m)}_{p,k}$ to $U^{(m)}_{p,k+1}$,  where the second index $k$ is modulo $m$.
In particular, the  preimage of $C_a$ has $m$ irreducible components (cuts)  for both $a=1,2$.
Denote these  preimages of  $C_1$ and $C_2$ by $C^{(m)}_{1, 1}, \dots, C^{(m)}_{1, m} \subset \tX^{(m)}$ and $C^{(m)}_{2, 1}, \dots, C^{(m)}_{2, m} \subset \tX^{(m)}$, respectively. 
Similarly, denote their projections  to $(X^{(m)}, 0)$ by $(D^{(m)}_{1, 1}, 0), \dots, (D^{(m)}_{1, m}, 0) \subset (X^{(m)}, 0)$ and $(D^{(m)}_{2, 1}, 0), \dots, (D^{(m)}_{2, m}, 0) \subset (X^{(m)}, 0)$. We know that the germs $(D^{(m)}_{1, k}, 0)$ and $(D^{(m)}_{2, k}, 0)$ are smooth for $1 \leq k \leq m$, since $(Z_{max}^{(m)})_{v_j} = 1$ for $1 \leq j \leq k_v$ (cf. (\ref{zmaxm})).

In the indexing  we assume that $C^{(m)}_{1,k}$ and  $C^{(m)}_{2,k}$ are subsets (cuts) in the same $U^{(m)}_{p,k}$ and  $D^{(m)}_{1,k}$ and  $D^{(m)}_{2,k}$
are  parametrized simultaneously from  $U^{(m)}_{p,k}$ via $\pi\circ \widetilde{f}_m$ (for the expressions see below).

  Note that $\mathbb{Z}_m$ acts in a cyclic way on the exceptional divisors $\{E^{(m)}_{v_j}\}_j$ as well.  Let us index them compatibly
such that the action  sends   $E^{(m)}_{v_j}$
 to  $E^{(m)}_{v_{j+1}}$, where the $j$ index is considered modulo $k_v$. Hence,  for any $k$ the curvettes  $C^{(m)}_{1,k+n_1k_v}$ and $  C^{(m)}_{1,k+n_2k_v}$ intersect the same exceptional divisor, for $n_1, \, n_2\in \mathbb{Z}$. (Note that $k_v|m$.)

If the simultaneous parametrization of $(D_1,0)$ and $(D_2,0)$ has the form (\ref{eq:para}) (where $C_1$ is given by $\beta=0$ and $C_2$ by $\beta=\beta_0\not=0$ with
$|\beta_0|\ll 1$), then in the covering above $U_p$  we have the equations $z^m=L=z_1=t^m=\alpha^m$. Hence $D^{(m)}_{1,k}$ and  $D^{(m)}_{2,k}$  have simultaneous parametrizations  from $U^{(m)}_{p,k}$ by
$$z=t\omega^k, \ z_1=t^m, \ \ z_i =\sum_{j\geq m}a_{i,j}(\beta)t^j \ (i\geq 2).$$
This, via $t\mapsto t\omega ^k$ transforms into
the following simultaneous parametrization of  $D^{(m)}_{1,k}$ and  $D^{(m)}_{2,k}$
\begin{equation}\label{eq:pauj}z=t, \ z_1=t^m, \ \ z_i =\sum_{j\geq m}a_{i,j}(\beta)\omega^{kj}t^j \ (i\geq 2).\end{equation}
Now, we can compare $q_v=q_0(D_1,D_2)$ and $q_{v_1}(D^{(m)}_{1,1},D^{(m)}_{2,1})$ via the definiton (\ref{eq:q12}). Recall that the multiplicity of both
$D^{(m)}_{1,1}$ and $D^{(m)}_{2,1}$ is 1, and, by the action, $q_0((D^{(m)}_{1,k},D^{(m)}_{2,l})= q_0(D^{(m)}_{1,k+1},D^{(m)}_{2,l+1})$.
Therefore, the discussion from subsection \ref{ss:simult} and comparison of (\ref{eq:Jl}) and (\ref{eq:pauj})  imply
\begin{equation*}\begin{split}
q_0(D_1,D_2)&=\frac{1}{m}\max_l q_0(D^{(m)}_{1,1},D^{(m)}_{2,l}) = \frac{1}{m}\max_lJ_l    \\
&=\frac{1}{m}J_m=\frac{1}{m}q_0(D^{(m)}_{1,m},D^{(m)}_{2,m})=\frac{1}{m}q^{(m)}_{v_m}=\frac{1}{m}q^{(m)}_{v_1}.\end{split}\end{equation*}
\end{proof}

Recall that the maximal ideal cycle $Z^{(m)}_{max}$ of $\tX^{(m)}$ has multiplicity one along each  $E^{(m)}_{v_i}$.
 Hence, in this way, 
for the singularity $(X^{(m)},0)$  and its resolution $\tX^{(m)}\to X^{(m)}$ we can apply Theorem \ref{innerrateone}. 
 Hence,  $q_{v_i}^{(m)} $ ---  in particular, $q_v$ too via Theorem \ref{ciclicinner} ---, can be recovered from the 
analytic semigroup of  $\tX^{(m)}$. 
However, this fact does not necessarily imply that they can be recovered from the semigroup of $\tX$ too!
In section \ref{ss:real}  we treat a case when such a reduction is possible. On the other hand, the example from section 5  shows that 
this is not possible in general.

\subsection{The $(Z_{max}, E_v)=0$ and  $m:=Z_{max,v}>1$ case with `unbranched vertex'  $v$}\label{ss:real} \

We recall our situation. We fix a  resolution $\pi:\tX\to X$ such that the line bundle $\calO_{\tX}(-Z_{max})$ has no basepoint and we also fix a vertex $v\in\mathcal{V}$ such that $m:=Z_{max,v}>1$. If $(Z_{max}, E_v)<0$ then $q_v=1$, cf.  (\ref{eq:q1}), hence in the sequel we assume that $(Z_{max}, E_v)=0$.
Write  $Z_{max}=\sum_ul_u E_u$. 

For any $u\in\mathcal{V}$ we denote by $\mathcal{V}(u)\subset \mathcal{V}$ the set  consisting  of $u$ and all its adjacent vertices in $\mathcal{T}$.

We consider the $m$ cyclic covering $f_m:X^{(m)}\to X$ branched along $\{L=0\}$, and the normalized pullback $\tilde{f}_m:\tX^{(m)}\to \tX$. This induces  a graph covering 
$\tau_m:\mathcal{T}^{(m)}\to \mathcal{T}$ of the corresponding dual graphs. 
Then, in general,  for any $w\in\mathcal{V}$,  $|\tau_m^{-1}(\{w\})|={\rm gcd}\{\{l_u\}_{u\in \mathcal{V}(w)}\}$, cf. \cite{cyclic,NBook}. In this subsection we assume that {\it  $\tau_m^{-1}(\{v\})$ has cardinality one}, that is, above the exceptional divisor $E_v$ of $\tX$ we have only one irreducible exceptional divisor in $\tX^{(m)}$, denoted by $E^{(m)}_{v_1}$. 

\smallskip

Our goal is to compute $q_{v_1}^{(m)}$ (hence, via Theorem \ref{ciclicinner},  $q_v$ too)
 in terms of the analytical semigroup $S_{an}$ of the resolution $\tX$ (after we eventually  modify $\tX$ by some additional blow ups, see below). 
The main ingredient is Theorem  \ref{innerrateone}, applied to  a   resolution of $X^{(m)}$. Recall  that $\tX^{(m)}$ is not a resolution of $X^{(m)}$, 
however it is smooth in a neighborhood of $E^{(m)}_{v_1}$ above which $\widehat{X}^{(m)}\to \tX^{(m)}$ is an isomorphism of smooth surfaces.
In particular, 
 $q_{v_1}^{(m)}$ can be computed already  in $\tX^{(m)}$. Moreover, 
 intersection with  $E^{(m)}_{v_1}$ of divisors supported on $E{(m)}=\pi_m^{-1}(0)$ make sense. 

From our assumption $(Z_{max}, E_v)=0$ we get that no strict transform of $L$ intersects $E_v$, hence, via the covering $\{L(x)=z^m\}$
no strict transform of $z$  intersects $E^{(m)}_{v_1}$.  
In particular, the inner rate of $E^{(m)}_{v_1}$ can be computed, via Theorem  \ref{innerrateone} applied for $E^{(m)}_{v_1}\subset \widehat{X}^{(m)}$,
 and reinterpreted in a neighbourhood of  $E^{(m)}_{v_1}\subset \tX^{(m)}$:
\begin{equation}\label{eq:coverq}
q_{v_1}^{(m)}=\min\,\{
l\,:\, \exists \ g^{(m)}\in \mathcal{O}_{(X^{(m)},0)},
 \ {\rm div}_{E^{(m)}_{v_1}}(\pi_m^*g^{(m)})=l, \  \  {\rm St}_{\tX^{(m)}}(g^{(m)}) \cap {E^{(m)\circ}_{v_1} \neq \emptyset} \}.
 \end{equation}
In the right hand side expression we would like to replace the condition  ${\rm St}_{\tX^{(m)}}(g^{(m)}) \cap {E^{(m)\circ}_{v_1} \neq \emptyset}$
 with 
a condition which is  readable from the analytic semigroup as well, namely with  a condition of type 
 $\ ( {\rm div}_{E^{(m)}}(\pi_m^*g^{(m)}), {E^{(m)}_{v_1}})<0$ (which make sense by the above observation).  

\bekezdes\label{bek:modif} First, however,  we need 
to  modify  the resolution just as in   subsection \ref{ss:AS}. 

We will modify the resolution $\tX$ as follows. Choose a function $g\in \mathcal{O}_{(X,0)}$ such that 
${\rm St}_{\tX}(g)\cap E_v^\circ \not=\emptyset$ (cf. Lemma \ref{lem:free} and  paragraph after it), and set 
$N_v:= {\rm div}_E(\pi^*g)_v$.  Fix also $r_v\in\mathbb{Q}_{>0}$ which satisfies Lemma
\ref{lem:iv} for $\pi$ and $v$,  and set an integer $b_v>r_v\cdot N_v$. 
Then we modify the resolution $\pi$ by blowing up each $\{E_v\cap E_u\}_{u\not=v}$ $mb_v$ times. (Observe the new  multiplicative factor $m$.)

For simplicity, this modified resolution will be denoted by the same symbol $\pi$ and $\tX$.

By the computation from subsection \ref{ss:cover} we obtain that  the vanishing order of $\tilde{f}_m^*\pi^*g=\pi_m^* f^*_mg$ along $E_{v_1}^{(m)}$ 
equals $N_v= {\rm div}_E(\pi^*g)_v$ too.   Note also that  $f_m^*g$ is one of the test functions in the right hand side of (\ref{eq:coverq}).
Hence,  for all function $g^{(m)}$ appearing in the right hand side of (\ref{eq:coverq}),  we can assume  that  ${\rm div}_{E^{(m)}_{v_1}}(\pi_m^*g^{(m)})=l\leq N_v$.

The analogue of Lemma \ref{lem:divst} is the following.

\begin{lemma}\label{lem:divst2}
If for some function $g^{(m)}$ of  $\mathcal{O}_{(X^{(m)},0)}$ the conditions ${\rm div}_{E^{(m)}_{v_1}}(\pi_m^*g^{(m)})\leq N_v$
 and
$ {\rm div}_{E^{(m)}}(\pi_m^*g^{(m)}), {E^{(m)}_{v_1}})<0$
are satisfied then 
 ${\rm St}_{\tX^{(m)}}(g^{(m)}) \cap {E^{(m)\circ}_{v_1} \neq \emptyset}$ as well. 
\end{lemma}
\begin{proof} Take a function  $g^{(m)}$ which satisfies the assumptions and 
write  ${\rm div}_{E^{(m)}_{v_1}}(\pi_m^*g^{(m)})=l$. 
Assume additionally that   ${\rm St}_{\tX^{(m)}}(g^{(m)}) \cap {E^{(m)\circ}_{v_1} =\emptyset}$.
Then 
 the strict transform ${\rm St}_{\tX^{(m)}}(g^{(m)})$ must intersect $E_{v_1}^{(m)}$ 
in an intersection point with some of its neighbours. Let $M$ be the $\Z_m$ action and consider the function 
$\prod_{i=0}^{m-1}  M^i g^{(m)}$. Since it is invariant, it equals $f_m ^*g$ for some $g\in \mathcal{O}_{(X,0)}$. Since the action preserves the vanishing orders, 
 ${\rm div}_{E^{(m)}_{v_1}}(\pi_m^*\prod _i M^ig^{(m)})=ml$, hence $ {\rm div}_{E}(\pi^*g)_v= ml\leq mN_v$ too. 
Furthermore,   ${\rm St}_{\tX}(g)$ passes at some  $E_v\cap E_u$ and 
in fact  $({\rm St}_{\tX}(g), E_v)\geq mb_v$ (since we performed $mb_v$ blow up).
 But then 
 $mN_v\geq  {\rm div}_{E}(\pi^*g)_v\geq ({\rm St}_{\tX}(g), E_v)/r_v \geq m b_v/r_v>m N_v$, a contradiction. 
\end{proof}

In conclusion, (\ref{eq:coverq}) can be rewritten as 
\begin{equation}\label{eq:coverq2}
q_{v_1}^{(m)}=\min\,\{
l\,:\, \exists \ g^{(m)}\in \mathcal{O}_{(X^{(m)},0)},
 \ {\rm div}_{E^{(m)}_{v_1}}(\pi_m^*g^{(m)})=l, \ ( {\rm div}_{E^{(m)}}(\pi_m^*g^{(m)}), {E^{(m)}_{v_1}})<0\}.
 \end{equation}

\bekezdes\label{bek:b1} 
 We wish to express $q_{v_1}^{(m)}$  in terms of the analytic semigroup associated with  the resolution $\tX$.

Below  it is convenient to identify $\mathcal{O}_{(X,0)}$ with the
 subring $f_m^*\mathcal{O}_{(X,0)}$ of $\mathcal{O}_{(X^{(m)},0)}$.

We pick a function   $g^{(m)}\in \mathcal{O}_{(X^{(m)},0)}$ with
  ${\rm div}_{E^{(m)}_{v_1}}(\pi_m^*g^{(m)})=l $ and $ ( {\rm div}_{E^{(m)}}(\pi_m^*g^{(m)}), {E^{(m)}_{v_1}})<0$. We say in such a case that $l$ is `realized by $g^{(m)}$'. We can also assume that $l\leq N_v$.

Since $X^{(m)}$ is a fiber product, $g^{(m)}$ can be written as $g^{(m)}=a_0+a_1z+\cdots a_{k-1}z^{m-1}$, with $a_k\in \mathcal{O}_{(X,0)}$. Let $D$ be the divisor of $\pi_m^*(g^{(m)})$, hence $D_{v_1}=l$.
Since we are searching for a minimal $l$, we can assume that $a_0\not=0$, otherwise we divide $g^{(m)}$ by a convenient $z$-power.

Let $M$ be the $\Z_m$ action  $(x,z)\mapsto (x, \omega z)$ acting on the local ring $\mathcal{O}_{(X^{(m)},0)}$
 and on the exceptional divisors of $\tX^{(m)}$ as well.
Then $M^ig^{(m)}=\sum_k a_k\omega ^{ki} z^k$ ($0\leq i<m$), whose divisor is $M^iD$. Then $D^*:=\min _{0\leq i<m}\{M^iD\}$ (the `symmetrized' $D$)
is the divisor of a {\it generic} combination $\sum _i c_i M^i g^{(m)}$ ($c_i\in\mathbb{C}$)
of the functions $M^ig^{(m)}$. But note that $$\sum _i c_i M^i g^{(m)}= \sum_i c_i \sum_k \omega^{ik}a_kz^k=
\sum_k(\sum_i c_i\omega^{ik})a_kz^k=\sum _k c_k' a_kz^k.$$
Since the matrix $(\omega^{ik})_{i,k}$ is invertible, the coefficients $c_k'=\sum_ic_i\omega^{ik}$ are generic, hence the new function
$g^{(m),s}=\sum _i c_i M^i g^{(m)}$ is a  generic combination of the functions $a_kz^k$ with divisor $D^*$.

The point is that the $\Z_m$ action  keeps the irreducible exceptional component  $E^{(m)}_{v_1}$ fixed
(since $k_v=1$), hence $(M^i D)_{v_1}=l$ for any $i$, hence $D^*_{v_1}=l$ too. Moreover,
$(D^*, E_{v_1}^{(m)})=(D, E_{v_1}^{(m)})+ (D^*-D, E_{v_1}^{(m)})<0$, hence the integer $l$ can be realized by this function $g^{(m),s}$ of type
$\sum_{k=0}^{m-1}c_k a_kz^k$ (and divisor $D^*$),  which, additionally, is  {\em  a generic combination of the summands
$a_kz^k$}. In particular,
\begin{equation}\label{divjo}
{\rm div}_{E^{(m)}_{v_1}} (a_kz^k)\geq l \ \ \mbox{for any $0\leq k<m$}.
\end{equation}

\begin{remark}\label{rem:key}
Before we continue the reduction of (\ref{eq:coverq2}), let us make a key observation. We would like to highlight the point  that distinguishes the 
case $|\tau_m^{-1}(\{v\})|=1$  from the general one. 
The crucial difference is given in the statements of the previous paragraphs. Namely, in the  $|\tau_m^{-1}(\{v\})|=1$  case
any testing function $g^{(m)}\in \mathcal{O}_{(X^{(m)},0)}$ can be replaced by a {\it generic combination} of germs of type $a_kz^k$. In fact, as we will see in the continuation, 
$g^{(m)}$ can be replaced even by a germ of type $a+cz^k$, where $a\in \mathcal{O}_{(X,0)}$ and $c$ is generic. However, in the 
$|\tau_m^{-1}(\{v\})|>1$ case, some special combination $a+cz^k$ might provide different divisor than the generic one (exactly this happens in the example from section 5).

In short, the existence of the above `monodromy symmetrization' makes the difference. 
\end{remark}
We divide the discussion into two cases.

\noindent {\bf Case I.} \ $a_k=0$ for all $k\geq 1$ and $a_0\not=0$.

\vspace{1mm}

\noindent {\bf Case II.} \ $a_k\not=0$ for at least one $k\geq 1$, and $a_0\not=0$.

In Case II 
we simplify  $g^{(m),s}=\sum_{k=0}^{m-1}c_k a_kz^k$ even further.
We set $D_k:={\rm div}_{E^{(m)}}( \pi_m^*a_kz^k)$, hence $\min_{k}D^k=D^*$.
Let $\mathcal{V}(v_1)$ be the set of vertices  consisting of $v_1$ and all its neighbors in $\mathcal{T}^{(m)}$.
The verification of the statements below are natural divisor manipulations and are left to the reader.

$\bullet$ 
Assume that $D^{k_0}_{v_1}>l$ for some $k_0\geq 1$. Then (the generic combination)
$$g^{(m),-}:=\sum _{0\leq k\leq m-1, k\not= k_0}c_k a_kz^k +c_{k_0}a_{k_0}z^{k_0-1}$$ has vanishing order $l$ along
$E^{(m)}_{v_1}$ (here we use the fact that the vanishing order of $z$  along $E_{v_1}^{(m)}$ is 1). Moreover, 
$( {\rm div}_{E^{(m)}}(\pi_m^*g^{(m),-}), {E^{(m)}_{v_1}}) \leq
( {\rm div}_{E^{(m)}}(\pi_m^*g^{(m),s}), {E^{(m)}_{v_1}})<0$. 
Hence, by repeating such steps we can find another function of type $\sum_k c_k a_kz^k$, a generic combination of summands $a_k z^k$,
which realizes $l$, and also $D^k_{v_1}=l$ for any $k\geq 1$ (whenever $a_kz^k$ is not zero).


$\bullet$
 If there exist two indices $1\leq k_1<k_2$ such that $a_{k_1}\not=0$, $a_{k_2}\not=0$
and $D^{k_1}|_{\mathcal{V}(v_1)}\not=D^{k_2}|_{\mathcal{V}(v_1)}$,  then a generic combination
$c_1a_{k_1}z^{k_1}+c_2a_{k_2}z^{k_2}$ also realizes $l$
(hence $c_1 a_{k_1}+c_2a_{k_2}z^{k_2-k_1}$ realizes $l-k_1$).

$\bullet$
 If all  $D^{k}|_{\mathcal{V}(v_1)}$ (associated with nonzero $a_k$'s) are all equal for all $k\geq 1$ then
$c_0a_0+c_ka_kz^k$ realizes $l$ (for some $k\geq 1$ with $a_k\not=0$).

$\bullet$
Hence, we can assume that if   $l$ is realized by some $g^{(m)}$, then $l$ (or a smaller value)  can be  realized  by a generic combination of
 type  $a+ c \cdot bz^k$, $a,b\in\mathcal{O}_{(X,0)}$, $c\in\mathbb{C}^*$. 
 By the realizability condition and by the fact that $a+ c \cdot bz^k$ is a generic combination, 
 one sees that the restriction of $z$ to
 $\mathcal{V}(v_1)$ cannot be a constant rational multiple of the restrictions to $\mathcal{V}(v_1)$ of both $a$ and $b$. 
If ${\rm div}_{E^{(m)}}(a)|_{\mathcal{V}(v_1)}$  
is not proportional to 
 ${\rm div}_{E^{(m)}}(z)|_{\mathcal{V}(v_1)}$  
then $a+c z^{k+k'}$ with $k':= {\rm div}_{E^{(m)}}(b)_{v_1}$
realizes the same value as   $a+ c \cdot bz^k$. Note that in this case $k+k'$ might be larger than $m$.
Otherwise,  $b+c z^{k'}$ 
realizes a smaller  value than    $a+ c \cdot bz^k$.

Next, we unify Cases I and II. 

If in Case I,  $a\in\mathcal{O}_{(X,0)}$ realizes some value $l$, then for some $c\in \mathbb{C}^*$ generic, $a+cz^l$ 
also realizes $l$.


\vspace{1mm}

  In particular, we can assume that in (\ref{eq:coverq2}) the minimal $l$ is realized by some function $g^{(m)}$ of type
 $a+cz^k$, where $a\in \mathcal{O}_{(X,0)}$, ${\rm div}_E(\pi^*a)_v={\rm div}_{E^{(m)}_{v_1}}(\pi_m^* a) \leq N_v$,   $k\geq 1$, 
and   $c\in\mathbb{C}^*$ is generic. 
For any fixed $a\in\mathcal{O}_{(X,0)}$  which appears in such a combination, 
there is an optimal $k$ which provides minimal $l$ (among functions of type $\{a+cz^k\}_k$, $a\in \mathcal{O}_{(X,0)}$, $c\in\mathbb{C}^*$):
  $$\min\{k\,:\, k=k\cdot {\rm div}_{E^{(m)}_{v_1}}(\pi_m^*z)\leq {\rm div}_{E^{(m)}_{v_1}}(\pi_m^*a), \
  k\cdot {\rm div}_{E^{(m)}}(\pi_m^*z)|_{\mathcal{V}(v_1)}\not\leq  {\rm div}_{E^{(m)}}(\pi_m^*a)|_{\mathcal{V}(v_1)}\}.$$
Here, the additional  inequality  $k\leq {\rm div}_{E^{(m)}_{v_1}}(\pi_m^*a)$  can be imposed by the following argument: if for some 
 $k> {\rm div}_{E^{(m)}_{v_1}}(\pi_m^*a)$  one has 
$  k\cdot {\rm div}_{E^{(m)}}(\pi_m^*z)|_{\mathcal{V}(v_1)}\not\leq  {\rm div}_{E^{(m)}}(\pi_m^*a)|_{\mathcal{V}(v_1)}$  then a similar inequality holds 
 for $k= {\rm div}_{E^{(m)}_{v_1}}(\pi_m^*a)$ as well. 

We summarize the output of our discussion in the following theorem.

  \begin{theorem}\label{th:4.2.6}  Assume that  $(Z_{max}, E_v)=0$ and  $m:=Z_{max,v}>1$.
Consider the resolution $\pi:\tX\to X$ as in \ref{bek:modif} and covering $\tX^{(m)}\to \tX$.
Then $$q_v=q^{(m)}_{v_1}/m, \ \mbox{ where} \ \  
q^{(m)}_{v_1}= \min\{ k(a)\ :\  a\in \mathcal{O}_{(X,0)}, \  {\rm div}_E(\pi^*a)_v\leq N_v\}, $$
where, for any $  a\in \mathcal{O}_{(X,0)}$, 
$$k(a):= 
 \min\{k\,:\, 1\leq k\leq {\rm div}_{E}(\pi^*a)_v, \
  k\cdot {\rm div}_{E^{(m)}}(\pi_m^*z)|_{\mathcal{V}(v_1)}\not\leq  {\rm div}_{E^{(m)}}(\pi_m^*a)|_{\mathcal{V}(v_1)}\}.$$
\end{theorem}

Next we wish to rewrite the statement of Theorem \ref{th:4.2.6} in terms of $\mathcal{S}_{an}(\tX)$,  where the resolution $\tX$
is as in \ref{bek:modif}, where all the intersection points $E_v\cap E_u$ are blown up sufficiently many times. 

Note that $z^m=L$, hence   ${\rm div}_{E^{(m)}}(\pi_m^*z)|_{\mathcal{V}(v_1)}=
\frac{1}{m}\cdot  {\rm div}_{E^{(m)}}(\pi_m^*L)|_{\mathcal{V}(v_1)}$. 
Moreover, the vanishing multiplicities 
$\{ {\rm div}_{E^{(m)}}(\pi_m^*a)\}_{u_i\in \tau_m^{-1}(u)}$  can be obtained from
$\{ {\rm div}_{E}(\pi^*a)\}_{u}$
via  multiplication by $m/(l_u,m)$, which is independent of $a$. 
Hence we obtain our final version:

  \begin{theorem}\label{th:4.2.7}  Assume that  $(Z_{max}, E_v)=0$ and  $m:=Z_{max,v}>1$.
Consider the resolution $\pi:\tX\to X$ as in \ref{bek:modif}. 
Then $m\cdot q_v=   \min\{ k(s)\ :\  s\in \mathcal{S}_{an}(\tX), \  s_v\leq N_v\}$, 
where
$$k(s):= 
 \min\{k\,:\, 1\leq k\leq s_v, \
  k\cdot Z_{max}|_{\mathcal{V}(v)}\not\leq m\cdot  s|_{\mathcal{V}(v)}\} \ \ \mbox{for any $  s\in \mathcal{S}_{an}(\tX)$}.$$
\end{theorem}




\section{The case of rational  singularities}

\subsection{}\label{ss:5.1} If a normal surface singularity is rational, then $\calS_{an}=\calS$ (see e.g. \cite{Lipman} or
\cite[Corolary 7.1.13]{NBook}), hence the analytic semigroup is determined by the topological type, i.e by
its resolution graph $\mathcal{T}$. Therefore,  it follows from the previous sections that whenever $(X,0)$ is rational and

$\bullet$ either $(Z_{min}, E_v)<0$,

$\bullet$ or  $(Z_{min}, E_v)=0$,  $Z_{min,v}=1$,

$\bullet$ or  $(Z_{min}, E_v)=0$,  $Z_{min,v}>1$, and $v$ is unibranched (i.e.   ${\rm gcd}\{ \{Z_{min, u}\}_{u \in \mathcal{V}(v)} \} = 1$)

\noindent then  the inner rate $q_v$ is topological.

In this section  we  show, that,  even in the rational singularity case, for branched vertices
$v \in \calv$ (i.e.  when
  ${\rm gcd}\{ \{Z_{min, u}\}_{u \in \mathcal{V}(v)} \} >1$),
 it may happen that the inner rate $q_v$ varies  when we modify the analytic type within a fixed topological type. 

Our example also shows (the rather surprising fact)  that for a rational singularity $(X,0)$, if we take  an $m$-branched covering $(X^{(m)},0)$ of $(X,0)$, branched along 
the generic linear section of $(X,0)$, the analytic semigroup of $X^{(m)},0)$ is not necessarily computable from the topology of $(X,0)$ and $m$.

Consider the rational resolution graph $\mathcal{T}$ from Figure 1 and  the vertex $v = v_{7}$ (for the meaning of the decorations see below). 
 We claim that there exist  two different  analytic structures
$\tX_1$, $\tX_2$  on the plumbed 4-manifold associated with the resolution graph $\mathcal{T}$, such that $q_{v, 1} = \frac{5}{3}$ and $q_{v, 2} \geq 2$.
Clearly, by contracting the exceptional divisors, we obtain two normal surface singularities $(X_1,0)$ and $(X_2,0)$  with the same topology.
The main difference between the analytic types will be the
position  in $E_{v_8}$  of the intersection points of $E_{v_8}$ with its adjacent curves $ E_{v_{1}}, \  E_{v_{2}}, \ E_{v_{3}},  \ E_{v_{4}}$.
 For $\tX_1$
we will choose a generic configuration of the $4$ intersection points of $E_{v_8}$, and for $\tX_2$ a special one, which will be described below.

\subsection{Some general statements regarding \underline{any} singularity $(X,0)$ with graph $\mathcal{T}$} \

Let $(X, 0) \subset \bC^N$ be the singularity with resolution $\tX$ and dual graph $\mathcal{T}$.
 We  fix  a generic linear function $L: (X, 0) \to \bC$, then ${\rm div}_E(\pi^*L) = Z_{max} =Z_{min} $, 
 which  is independent of
 the analytic type, cf. \cite{Artin62,Artin66}.
 Figure 1 shows the graph $\mathcal{T}$ of $\tX$,
the integers in $[\, \ ]$ are the self-intersection numbers, the pair of integers in $(\ )$ are the multiplicities of $Z_{min}$  and of ${\rm div}_E(\pi^*y)$, where $y$ is a function which will be  defined below.

For any  $u\in \mathcal{V}(v)$  one has
 $(Z_{min}, E_{u}) = 0$ and $\pi^*L \not\in  H^0(\calO_{\tX}(-Z_{min} - E_u))$, hence
from the exact sequence $0 \to H^0(\calO_{\tX}(-Z_{min} - E_u)) \to H^0(\calO_{\tX}(-Z_{min})) \to H^0(\calO_{E_u}(-Z_{min}))\simeq  \bC$
we obtain that
$$\dim\left(H^0(\calO_{\tX}(-Z_{min}))/H^0(\calO_{\tX}(-Z_{min} - E_u))\right) = 1.$$
 Then there exists a unique minimal cycle $Z_1\in\calS$ such that  $Z_1\geq Z_{min}+E_u$ ($u\in \mathcal{V}(v)$). In fact,
by Generalized Laufer Algorithm (cf. \ref{ss:notation}) the cycle $Z_1$ is the same for any $u\in \mathcal{V}(v)$ (and
 $Z_1$ is also the minimal element in $\calS$, such that $Z_1 \geq  4E_v$).
(In fact, the algorithm shows that $s(Z_1+E_u)=s(Z_1+E_v)=s(4E_v)$. This is the cycle $Z_1$.  For the $s$-operator see subsection \ref{ss:notation}.)

Since $(X, 0)$ is rational, thus $\calS = \calS_{an}$, $Z_1={\rm div}_E(\pi^*y)$ for some function $y\in\mathcal{O}_{(X,0)}$.
The identity $\calS=\calS_{an}$ implies that $H^0(\calO_{\tX}(-Z_{min} - E_u)) = H^0(\calO_{\tX}(- Z_1))$ for any $u\in \mathcal{V}(v)$.
From Figure 1 we also read the multiplicities
$(Z_1)_u$ for  $u\in \mathcal{V}(v) = \{v=v_7, v_{1}, v_{5}, v_{6}\}$, they are   4,8,4,4 respectively. This shows that the functions $\pi^*L$ and $\pi^*y$ are 
proportional  on $\mathcal{V}(v)$:
 $(Z_{1})_{\mathcal{V}(v)} = \frac{4}{3} \cdot (Z_{min})_{ \mathcal{V}(v)}$.  Also,   $(Z_1,E_u)=0$,  similarly as $(Z_{min},E_u)=0$ for any $u\in \mathcal{V}(v)$.
From the above exact sequence and from
 $0 \to H^0(\calO_{\tX}(-Z_1 - E_u)) \to H^0(\calO_{\tX}(-Z_1)) \to H^0(\calO_{E_u}(-Z_1))\simeq  \bC$
we obtain  $\dim(H^0(\calO_{\tX}(-Z_1))/ H^0(\calO_{\tX}(-Z_1 - E_u)))=1$ generated by $\pi^*y$. But  $H^0(\calO_{\tX}(-Z_1 - E_u))=
 H^0(\calO_{\tX}(-s(Z_1 + E_u)))$. By Generalized Laufer Algorithm
$s(Z_1+E_u)=s(Z_1+E_v)=s(Z_1+\sum _{u\in \mathcal{V}(v)}E_u)$ (once some  $E_u$ is added to $Z_1$, in the next step we can  add $E_v$ as well, hence
in continuation we add  all the vertices of $\mathcal{V}(v)$ step by step.) Hence
\begin{equation}\label{eq:3}
H^0(\calO_{\tX})/H^0(\calO_{\tX}(-(Z_1+E)|_{\mathcal{V}(v)})) \ \ \mbox{is freely generated by the cosets of
$\{1, \pi^*L, \pi^* y\}$}.\end{equation}
Consider next  $(X^{(3)}, 0)$, the $3$-th cyclic cover of $(X, 0)$ ramified along $\{L=0\}$ and its partial resolution $\tX^{(3)}$
 with exceptional set $E^{(3)}$.
Recall  that  $(X^{(3)}, 0) \subset \{z^3 = L(x)\} \subset (\bC^{N+1},0)$,  and  ${\rm div}_{E^{(3)}}(\pi_3^*z) = Z_{max}^{(3)}$, cf. (\ref{zmaxm}).
 Denote also $Z^{(3)}_1 ={\rm  div}_{E^{(m)}}(\tilde{f}_3^*\pi^*y)$. Above $E_v$ there are three irreducible components of $E^{(3)}$: $E^{(3)}_{\tilde{v}_1}$,
$E^{(3)}_{\tilde{v}_{2}}$ and $E^{(3)}_{\tilde{v}_{3}}$. One verifies  (using the formulae from subsection \ref{ss:cover}) that
$Z^{(3)}_1|_{\mathcal{V}(\tilde{v}_1)}=4\cdot Z_{max}^{(3)}|_{\mathcal{V}(\tilde{v}_1)}$, for the concrete values see
Figure 3.  Also, $(Z_{max}^{(3)}, E^{(3)}_{\tilde{v}_1})=0$.

Recall  that $\calO_{(X^{(3)}, 0)}$ is generated by $f_3^*(\calO_{(X, 0)})$ and $z$, hence any function of  $\calO_{(X^{(3)}, 0)}$
has the form $a_0+za_1+a_2z^2$, where $a_i\in f_3^*(\calO_{(X, 0)})$. Since $({\rm div}_{E^{(3)}} \tilde{f}_3^*(\pi^*a))_{\tilde{v}_1} = ({\rm div}_E(\pi^*a))_v$ for
any $a\in \mathcal{O}_{(X,0)}$, (\ref{eq:3}) imply that
\begin{equation}\label{eq:3b}H^0(\calO_{\tX^{(3)}})/H^0(\calO_{\tX^{(3)}}(-(Z^{(3)}_1+E^{(3)})|_{\mathcal{V}({\tilde{v}_1} )})) \ \mbox{is generated by $(1, \pi_3^*z, \ldots,\pi_3^* z^4,
\tilde{f}_3^*\pi^*y)$.}\end{equation}
\begin{lemma}\label{lem:5}  For any analytic type supported on $\mathcal{T}$  one has
 $q_{\tilde{v}_{ 1}} \geq 5$.
\end{lemma}
\begin{proof}  Since $(Z_{max}^{(3)}, E^{(3)}_{\tilde{v}_1})=(Z_1^{(3)}, E^{(3)}_{\tilde{v}_1})=0$,  by (\ref{eq:3b}), for any
  $\tilde{g}\in \mathcal{O}_{X^{(3)}}$ such that ${\rm div}_{E^{(3)}}(\pi_3^*\tilde{g})_{\tilde{v}_1}\leq 4$ one also has
$({\rm div}_{E^{(3)}}(\pi_3^*\tilde{g}), E^{(3)}_{\tilde{v}_1})= 0$. (One should only test  $ \pi_3^*z, \ldots,\pi_3^* z^4,
\tilde{f}_3^*\pi^*y$.)
%
\end{proof}

Next,  we analyse the preimage of $E_{v_8}$.
Since ${\rm gcd} \{Z_{min,u}\,|\, u\in \mathcal{V}(v_8)\}={\rm gcd}\{6, 3, 2, 5\} = 1$, above $E_{v_8}$ there is only one irreducible exceptional divisor of $E^{(3)}$, denoted by
$E^{(3)}_{\tilde{v}_8}$. The covering $E^{(3)}_{\tilde{v}_8}\to E_{v_8}$ has degree three, and it has three branch points $p_2, p_3,p_4\in E_{v_8}$,
the intersection points with
$E_{v_2}, E_{v_3}, E_{v_{4}}$ respectively. An Euler characteristic computation  shows that $E^{(3)}_{\tilde{v}_8}$ is an elliptic curve.

We also write $p_1:= E_{v_8} \cap E_{v_{1}}\in E_{v_8}$ as well. (It is not a branch point since ${\rm gcd}\{Z_{min, v_8}, Z_{min, v_{1}}\}= {\rm gcd}\{3,6\}=3$, hence
the  preimage of $p_1$ consists of  three points.)



\subsection{Discussion regarding the elliptic exceptional divisor $E^{(3)}_{\tilde{v}_8}\subset E^{(3)}$}\label{ss:ell} \

Recall that  we wish to construct a special analytic type $(X_2, 0)$ for which $q_v \geq \frac{6}{3} = 2$. For this we  fix the intersection points $p_2, \
 p_3,\  p_4 $, however  we will choose the intersection point $p_1$ in a special way (with respect to $\{p_2,p_3,p_4\}$). This special position of $p_1$ will be determined  below.

Let us identify $E_{v_8}$ with $\mathbb{P}^1$,  and the points $p_3, \, p_4, \,  p_2$ with $0,\, 1,\, \infty$.   We claim that there exists a unique elliptic curve $C$ and covering
$C \to \mathbb{P}^1$ of degree 3 with branch points  $0,\, 1,\, \infty$.

Indeed,
  such  coverings correspond to  surjective monodromy representations $\rho: \pi_1(\mathbb{P}^1 \setminus \{0, 1, \infty\}, *) \to S_3$
(where $S_3$ is the permutation group of 3 elements)
up to global conjugation.
The fundamental group is generated by elementary  loops $\gamma_2, \gamma_3, \gamma_4$ (`around $p_2,p_3,p_4$')
subject to $\gamma_2 \gamma_3 \gamma_4 = 1$. The image of any $\gamma_i$ is a
 $3$-cycle, hence it equals $(123)^{a_i}$ for some $a_i\in\{1,2\}$.
Since  $\gamma_2 \gamma_3 \gamma_4 = 1$ we have $a_2+a_3+a_4\equiv 0 \ ({\rm mod} \, 3)$. Hence $(a_1,a_2,a_3)$ equals either $(1,1,1)$ or $(2,2,2)$. Since these two
possible choices give conjugate representations, the uniqueness is proved. Let $\tilde{0},\, \tilde{1}, \, \tilde{\infty}$ be the points in $C$ above $ 0,\, 1,\, \infty$
(they correspond in $E^{(3)}_{\tilde{v}_8}$ to the points  $\tilde{p}_3,\,  \tilde{p}_4,\,  \tilde{p}_2$).

The  covering can be realized by $C:=\{[a:b:c]\in \mathbb{P}^2\ |\ c^3=a(a-b)b\}\to \mathbb{P}^1$, $[a:b:c]\mapsto [a:b]$.
Set $\calL:= \mathcal{O}_C(2\tilde{\infty})$, where $\tilde{\infty}:= [1:0:0]$.  By Riemann-Roch theorem $h^0(C,\calL)=2$.
Any divisor linearly equivalent to $2\tilde{\infty}$  has the form $P+Q$, where $P+Q+\tilde{\infty}$ are the intersection points of $C$ with the line $\{\alpha c=\beta b\}_{[\alpha:\beta]\in\mathbb{P}^1}$. Let $s_1\in H^0(C,\calL)$ be the section corresponding to $\{b=0\}$ and  with divisor $2\tilde{\infty}$, while $s_2\in H^0(C,\calL)$ corresponds to
 $\{c=0\}$  with  ${\rm div}(s_2)=\tilde{0}+\tilde{1}$. Here $\tilde{0}=[0:1:0]$ and $\tilde{1}=[1:1:0]$. 
We wish to find those
sections for which $P=Q$. Clearly, one possibility is realized by $\{b=0\}$ when $P=Q=\tilde{\infty}$. Otherwise we take the affine coordinates
$c/b=\lambda$ and $a/b=r$, we  rewrite $\ c^3=a(a-b)b$ into the equation $\lambda^3=r(r-1)$ and we wish to find $\lambda $  such that the equation has a double 
root in $r$.
This happens exactly when $\lambda^3+1/4=0$ and  $r=1/2$. This means that if we consider the point $[a:b]=[1:2]\in \mathbb{P}^1$,
  and the three points $\tilde{P}_1,\, \tilde{P}_2, \ \tilde{P}_3$ above it, then for the three solutions of $\lambda^3+1/4=0$ the divisor of the corresponding section $s_2-\lambda s_1$
in $H^0(C,\calL)$ is $2\tilde{P}_i$.  For other values $\lambda\in \bC$,  ${\rm div}(s_2-\lambda s_1)$ consists of two distinct points.
Since in this case, $r(r-1)=\lambda^3$ has two distinct roots in $r=a/b$, the projection of these two  point  are distinct in $\mathbb{P}^1$ too.

(The fact that there are exactly four sections  with divisor of  the form $2\tilde{P}$ corresponds to the fact that in ${\rm Pic}(C)$ --- that is, in the 2-torus ---
there are exactly four torsion points of order two.)

\subsection{The construction of the different analytic types} \

The resolution space $\tX$ can be constructed by analytic plumbing corresponding to the resolution graph.
Since  $E_{v_8}\simeq {\mathbb P}^1$, any three points of it can be transformed by an automorphism of $E_{v_8}$ to $0,1,\infty$.
However, the choice of a fourth point creates a 1-dimensional moduli. Let us fix the intersection points
$p_2, \, p_3$ and $p_4$ of $E_{v_8}$. Then by moving the intersection point $p_1=E_{v_8} \cap E_{v_{1}}$  in the plumbing construction we can modify the analytic structure of $\tX$.

The  analytic type $(X_1, 0)$  will be realized by a generic choice of $p_1$ (with respect to $p_2,p_3,p_4$), and in this case we will prove that  $q_v = \frac{5}{3}$
(or, equivalently,  $q_{\tilde{v}_{1}} = 5$). The analytic type  $(X_2, 0)$  will be identified by the  special choice of
$p_1=[1:2]$  for which $q_v \geq \frac{6}{3} = 2$ (or, equivalently, $q_{\tilde{v}_{ 1}} \geq 6$).

\subsection{The generic choice of $p_1$, the singularity $(X_1,0)$}\label{ss:x1} \

 Let us start with the first case, when the choice of
$p_1$ is generic.
Denote the preimages of $p_2, p_3, p_4$ under the cyclic cover by $\tilde{p}_2, \tilde{p}_3, \tilde{p}_4 \in E^{(3)}_{\tilde{v}_{8}}$ and the preimages of $p_1$ by $\tilde{p}_{1, 1} = E^{(3)}_{\tilde{v}_{1, 1}} \cap E^{(3)}_{\tilde{v}_{8}}, \  \tilde{p}_{1, 2} = E^{(3)}_{\tilde{v}_{1, 2}} \cap E^{(3)}_{\tilde{v}_{8}}, \
 \tilde{p}_{1, 3} = E^{(3)}_{\tilde{v}_{1, 3}} \cap E^{(3)}_{\tilde{v}_{8}}$.
The divisors of the pullback of the functions $z^4$ and $  f_3^*y$ restricted to $E^{(3)}_{\tilde{v}_8}$
give the trivial line bundle $\calO_{E^{(3)}_{\tilde{v}_{8}}}$, which implies:
$$\calO_{E^{(3)}_{\tilde{v}_{8}}}\Big(4 E^{(3)}_{\tilde{v}_{8}} + 8 \sum_{1 \leq i \leq 3} \tilde{p}_{1, i} + 20\tilde{p}_2 + 8\tilde{p}_3 + 8  \tilde{p}_4 \Big) \cong
\calO_{E^{(3)}_{\tilde{v}_{8}}} \Big(4 E^{(3)}_{\tilde{v}_{8}} + 8 \sum_{1 \leq i \leq 3} \tilde{p}_{1, i} +  18 \tilde{p}_2 + 9 \tilde{p}_3 +  9 \tilde{p}_4\Big) = \calO_{E^{(3)}_{\tilde{v}_{8}}}.$$
In particular, $\calO_{E^{(3)}_{\tilde{v}_{8}}}(2 \tilde{p}_2) = \calO_{E^{(3)}_{\tilde{v}_{8}}}( \tilde{p}_3 + \tilde{p}_4)$; denote this line bundle by  $\calL$.

Consider a section $s_1 \in H^0(E^{(3)}_{\tilde{v}_{8}}, \calL)$ with divisor $ 2 \tilde{p}_2$ and $s_2 = s_1 \cdot \frac{\tilde{f_3}^*(y)}{z^4} \in H^0(E^{(3)}_{\tilde{v}_{8}}, \calL)$
with divisor  $ \tilde{p}_3 + \tilde{p}_4$.  They generate  the complete linear system $\langle s_1, s_2\rangle
 = H^0(E^{(3)}_{\tilde{v}_{8}}, \calL)$.  Since  $s_1, s_2$ do not vanish at $\tilde{p}_{1, 1}$, there  exists a   unique parameter $ 0 \neq \mu \in \bC$, such that $s = s_1 + \mu s_2$ vanishes at $\tilde{p}_{1, 1}$.
By the discussion from \ref{ss:ell}, the other zero of $s$, say $\tilde{p}^*$, is different from all the points $\{\tilde{p}_j\}_{j=2,3,4}$ and $\{\tilde{p}_{1,i}\}_{i=1,2,3}$.  
Moreover, the vanishing order of  $s$ at $\tilde{p}_{1, 1}$ is $1$
(cf.  the description \ref{ss:ell}).


Consider the function $g = z^4 + \mu f_3^*(y) = z^4 \cdot \frac{s}{s_1}\in \mathcal{O}_{X^{(3)}}$.
Since the vanishing order of $\pi_3^*z^4 $ (respectively of $\pi^*_3f_3^*y$)  along  the exceptional divisors $E^{(3)}_{\tilde{v}_{8}},\
 E^{(3)}_{\tilde{v}_{2}},\  E^{(3)}_{\tilde{v}_{3}}$ and $ E^{(3)}_{\tilde{v}_{4}}$ is $4, 20, 8,8$ (respectively $4, 18, 9,9$),
the vanishing order of  $\pi^*_3g$ along these divisors
is  $4, 18, 8, 8$ respectively.
Moreover, $g$ has a strict transform passing at $\tilde{p}^*$, and 
 the intersection multiplicity of $\frac{g}{z^4}$ with $E^{(3)}_{\tilde{v}_{8}}$ at $\tilde{p}_{1, 1}$ is exactly  $1$.

We claim that this intersection multiplicity 1 cannot be realized by a strict transform  of $g$ at $\tilde{p}_{1,1}$. Indeed, if such a strict transform existed, then the sum of
vanishing orders along the  adjacent vertices of $\tilde{v}_{1.1}$ (including this strict transform as well) would be at least $4+4+8+1=17$. Since the self-intersection of
$ E^{(3)}_{\tilde{v}_{1, 1}}$
 is $-2$,
we necessarily would have that the  vanishing order of $\pi_3^*g$ along   $ E^{(3)}_{\tilde{v}_{1, 1}}$ is $\geq 9$. But then
the intersection multiplicity of $\frac{g}{z^4}$ with $E^{(3)}_{\tilde{v}_{8}}$ at $\tilde{p}_{1, 1}$ would be  $\geq (9-8)+1= 2$.

Hence $g$ has no strict transform at $\tilde{p}_{1,1}$ and the vanishing order of the function $\pi_3^*g$ along the exceptional divisor $E^{(3)}_{\tilde{v}_{1, 1}}$ is exactly $9$.

\medskip

Next, we analyze the neighborhood of $E^{(3)}_{\tilde{v}_1}$. The vanishing orders of the adjacent curves is (at least)  $4+4+9=17$, the self-intersection is $-4$.
Therefore, $\pi^*_3g$ must vanish along the exceptional divisor $E^{(3)} _{\tilde{v}_{1}}$ of order
at least $5$. On the other hand,  the vanishing orders along  any of
$E^{(3)} _{\tilde{v}_{1}}$,  $E^{(3)} _{\tilde{v}_{5,1}}$  $E^{(3)} _{\tilde{v}_{6,1}}$ of $\pi_3^*z^5$ is 5, hence
 for a generic parameter $\gamma \in \bC$, the function $\tilde{g}:=\pi^*_3(g + \gamma \cdot z^5)$ vanishes along the exceptional divisor $E^{(3)}_{\tilde{v}_{ 1}}$ of
 order 5 and along its neighbors   $E^{(3)} _{\tilde{v}_{5,1}}$,  $E^{(3)} _{\tilde{v}_{6,1}}$  $E^{(3)} _{\tilde{v}_{1,1}}$ of order $\leq 5, \leq 5, 9$. In particular,
$({\rm div}_{E^{(3)}}(\pi^*_3(g + \gamma \cdot z^5)), E^{(3)}_{\tilde{v}_{ 1}}) < 0$.

We claim that in fact  ${\rm St}_{\tX^{(3)}}(\tilde{g}) \cap E^{(3)\circ}_{\tilde{v}_1} \neq \emptyset$. Indeed,  notice that the vanishing order of $\pi^*_3g$
along the exceptional divisor $E^{(3)} _{\tilde{v}_{1}}$ is at leats 5, which implies via Laufer algorithm that the vanishing orders along 
$  E^{(3)} _{\tilde{v}_{5,1}}, E^{(3)} _{\tilde{v}_{6,1}}$ are also  at least $5$. Furthermore,  $z^5$ has no strict transform going through
$E^{(3)} _{\tilde{v}_{1}} \cap E^{(3)} _{\tilde{v}_{5,1}}, E^{(3)} _{\tilde{v}_{1}} \cap E^{(3)} _{\tilde{v}_{6,1}}$. 
This means, that for a generic parameter $\gamma \in \bC$, the function $\tilde{g}$ has no strict transform passing through  the  points $E^{(3)} _{\tilde{v}_{1}} \cap E^{(3)} _{\tilde{v}_{5,1}}, E^{(3)} _{\tilde{v}_{1}} \cap E^{(3)} _{\tilde{v}_{6,1}}$. Hence, the only remaining  possibility would be that ${\rm St}_{\tX^{(3)}}(\tilde{g})$ goes through $E^{(3)} _{\tilde{v}_{1}} \cap E^{(3)} _{\tilde{v}_{1,1}}$. However,  the vanishing order of $\tilde{g}$ along $E^{(3)} _{\tilde{v}_{1,1}}$ is $9$ (which is a vertex with selfintersection $-2$) and along the exceptional divisors $E^{(3)}_{\tilde{v}_{ 7}}, E^{(3)}_{\tilde{v}_{ 8}}, E^{(3)}_{\tilde{v}_{ 9}}$ at least $5, 4, 9$, which excludes the possibility of the 
existence of such a strict transform.

Hence, by Theorem \ref{innerrateone} and Lemma \ref{lem:5},
 for the analytic type $(X_1, 0)$, one has $q_{\tilde{v}_{1}} = 5$, or equivalently $q_v = \frac{5}{3}$.

\subsection{The special choice of $p_1$, the singularity $(X_2,0)$}\

Let us  identify again   $E_{v_8}$ with $\mathbb{P}^1$,  and the points $p_3, \, p_4, \,  p_2$ with $0,\, 1,\, \infty$. In our second construction the
plumbing construction is performed  in such a way that $p_1=E_{v_8}\cap E_{v_{1}}$ is $1/2$, cf. \ref{ss:ell}. Then the points $\tilde{P}_i$ (considered in \ref{ss:ell})
correspond to the points
$\tilde{p}_{1,i}$ in the covering $\widetilde{X}^{(m)}$.

 Let us consider again the pencil    $g = z^4 + \mu \tilde{f}_3^*(y)$. It vanishes along
 $E^{(3)}_{\tilde{v}_{8}},\
 E^{(3)}_{\tilde{v}_{2}},\  E^{(3)}_{\tilde{v}_{3}}$ and $ E^{(3)}_{\tilde{v}_{4}}$ of order $4, 18, 8, 8$ respectively.
Moreover, by the discussion from subsection \ref{ss:ell},   for three special  values  $\mu_i$
the intersection multiplicity of $\pi_3^*(\frac{g}{z^4})$ with $E^{(3)}_{\tilde{v}_{8}}$ at $\tilde{p}_{1, i}$ equals $2$   ($1\leq i\leq 3$).

We fix $\mu_1$.  Hence, by this choice of the analytic plumbing (i.e. of the analytic structure $(X_2,0)$)
 the intersection multiplicity of $\pi^*_3(\frac{g}{z^4})$ with $E^{(3)}_{\tilde{v}_{8}}$ at $\tilde{p}_{1, 1}$ equals $2$.

\medskip

Similarly as in subsection \ref{ss:x1}, this fact implies that  the vanishing order of $\pi^*_3g$ along $E^{(3)}_{\tilde{v}_{1, 1}}$  must increase from 8:
it must be  at least $9$ and at most $10$,
 and if the vanishing order of $\pi_3^*g$ along $E^{(3)}_{\tilde{v}_{1, 1}}$ is $9$, then $\pi_3^*g$ must have a strict transform passing through the intersection point $E^{(3)}_{\tilde{v}_{8}} \cap E^{(3)}_{\tilde{v}_{1, 1}} = \tilde{p}_{1, 1}$.

\medskip

In this paragraph we exclude the second possibility, namely when the vanishing order of $\pi^*_3g$ along $E^{(3)}_{\tilde{v}_{1, 1}}$ is $9$ and it has a strict transform passing at
$\tilde{p}_{1,1}$. Indeed, in that case, (similarly as in subsection \ref{ss:x1}), via Laufer Algorithm,  the vansihing order of $\pi^*_3g$ along $E^{(3)}_{\tilde{v}_1}$ should be at least $5$ and along  $E^{(3)}_{\tilde{v}_{9,1}}$ should be at least $9$. Hence, counting the sum of multiplicities of the neighbours of
 $E^{(3)}_{\tilde{v}_{1,1}}$ (including the strict transform at $\tilde{p}_{1,1}$) we get  $\geq (4+5+9+1)=19$. Since the selfintersection and
the multiplicity  of  $E^{(3)}_{\tilde{v}_{1,1}}$   is $(-2)$  and 9, we obtain $({\rm div}(\pi^*_3g), E^{(3)}_{\tilde{v}_{1,1}})>0$, a contradiction.

\medskip

In particular, regarding the function $\pi^*_3g$,  we can assume that its vanishing order
 along  $E^{(3)}_{\tilde{v}_{1, 1}}$ is $10$. Next we focus on the exceptional divisor   $E^{(3)}_{\tilde{v}_{1}}$ (on which we wish to test  $q_{\tilde{v}_1}$).
First,  we show that a possible $q_{\tilde{v}_1}=5 $ cannot be realized by the function $g$, that is, simultaneously ${\rm div}(\pi^*_3g)_{\tilde{v}_{1}}=5$ and $({\rm div}_{E^{(3)}} (\pi^*_3g), E^{(3)}_{\tilde{v}_{1}})<0$   (cf. Theorem \ref{innerrateone})
cannot happen.

Indeed, since the vanishing order of $\pi^*_3g$ along  $E^{(3)}_{\tilde{v}_{1, 1}}$ is $10$, by the Generalized Laufer Algorithm
the vanishing order along  its neighbor  $E^{(3)}_{\tilde{v}_{1}}$ should be at least 5, and let us assume that indeed it is 5.
Then the same algorithm applied for the neighbors   $E^{(3)}_{\tilde{v}_{5, 1}}$ and    $E^{(3)}_{\tilde{v}_{6, 1}}$ of
$E^{(3)}_{\tilde{v}_{1}}$ gives that their vanishing orders should be at least 5. Hence, the sum of vanishing orders of the neighbors of
 $E^{(3)}_{\tilde{v}_{1}}$  is $\geq (5+5+10)=20$, its vanishing order is 5 and selfintersection $(-4)$, hence
 $({\rm div}_{E^{(3)}} (\pi^*_3g), E^{(3)}_{\tilde{v}_{1}})\geq 0$.

Finally, we wish to show that  $q_{\tilde{v}_1}=5 $ cannot be realized by any function, that is,  there is no function $h\in\mathcal{O}_{(X_2^{(3)},0)}$ such that
${\rm div}(\pi^*_3h)_{\tilde{v}_{1}}=5$ and   $({\rm div}_{E^{(3)}} (\pi^*_3h), E^{(3)}_{\tilde{v}_{1}})<0$.
This (together with Lemma \ref{lem:5}) would imply that $q_{\tilde{v}_{1}} \geq 6$.

\medskip

Assume that such a function  $h \in \calO_{(\tX^{(3)}_2,0)}$ exists.
We focus on its vanishing order at  $E^{(3)}_{\tilde{v}_{1}}$, which is 5, and the vanishing orders of its neighbors
  $E^{(3)} _{\tilde{v}_{5,1}}$  $E^{(3)} _{\tilde{v}_{6,1}}$  $E^{(3)} _{\tilde{v}_{1,1}}$. The vanishing orders of $\pi_3^*z^5 $ along them are 5,5,10, hence considering a generic combination of $h$ and $z^5$ we can assume that the vanishing order of $h$ along this three curves is not greater than 5,5,10. But if they are 5,5,10, then
 $({\rm div}_{E^{(3)}} (\pi^*_3h), E^{(3)}_{\tilde{v}_{1}})=0$, so at least one of them should be smaller.

Note also that if the vanishing order at  $E^{(3)}_{\tilde{v}_{1}}$ is  5, then by the Generalized Laufer Algorithm, the  vanishing orders of
  $E^{(3)} _{\tilde{v}_{5,1}}$  $E^{(3)} _{\tilde{v}_{6,1}}$  $E^{(3)} _{\tilde{v}_{1,1}}$ should be $\geq 5, \geq 5, \geq 9$.
Hence if such a function $h$ exists then its vanishing orders along the neighbors should be $5,5,9$. (We will show that `this 9' is not possible.)

Recall that $h$ has the form $\sum_kf_3^*(a_k)z^k$ where $a_k\in \mathcal{O}_{(X_2,0)}$.
In this sum there are terms for which we automatically know that
their vanishing orders along  $E^{(3)} _{\tilde{v}_{5,1}}$  $E^{(3)} _{\tilde{v}_{6,1}}$  $E^{(3)} _{\tilde{v}_{1,1}}$ are greater than 5,5,10.
Some of them can be seen immediately: type $z^k$ with $k\geq 5$, or $z^kf_3^*y$ with $k\geq 1$. Also, by Generalized Laufer Algorithm we see that if
some $a\in \mathcal{O}_{(X_2,0)}$  and  $\pi^*a\in H^0(\tX, \mathcal{O}_{\tX}(-5E_v))$, then it also vanish along $E_{v_5}, E_{v_6}, E_{v_1}$  with orders 5,5,10,  hence the same will be true for their lifts to the cyclic covering. 

Denote the sum of those terms by $h'$, then $h$ can be written as the sum of $h'$ with a linear combination of terms of type $ z,z^2,z^3,z^4, f_3^*y\in \mathcal{O}_{X_2^{(3)}}$.
The terms $\{z^k\}_{k=1,2,3} $ cannot appear since in that case the vanishing order of $h$ along $E_{\tilde{v}_{1}}$ would be $\leq 3$.  Hence $h = h' + c_1z^4+c_2f_3^*y$, where $c_1,c_2\in \bC$.

We already verified that for the special member of the pencil, namely for $g$, the vanishing order along $E^{(3)} _{\tilde{v}_{1,1}}$ is 10 (hence not 9), which gives that the vanishing order of $h$ along $E_{\tilde{v}_{1}}$ would be $10$.
On the other hand, by the construction from subsection \ref{ss:ell} we see that for any other member of the pencil, the
 vanishing order along   $E^{(3)} _{\tilde{v}_{1,1}}$ is $8$, since the divisor of the corresponding $s\in H^0(C,\calL)$ is not supported at $\tilde{p}_{1,1}$. In that case the vanishing order of $h$ along $E_{\tilde{v}_{1}}$ would be $8$, which is again a contradiction. This finishes the proof completely.

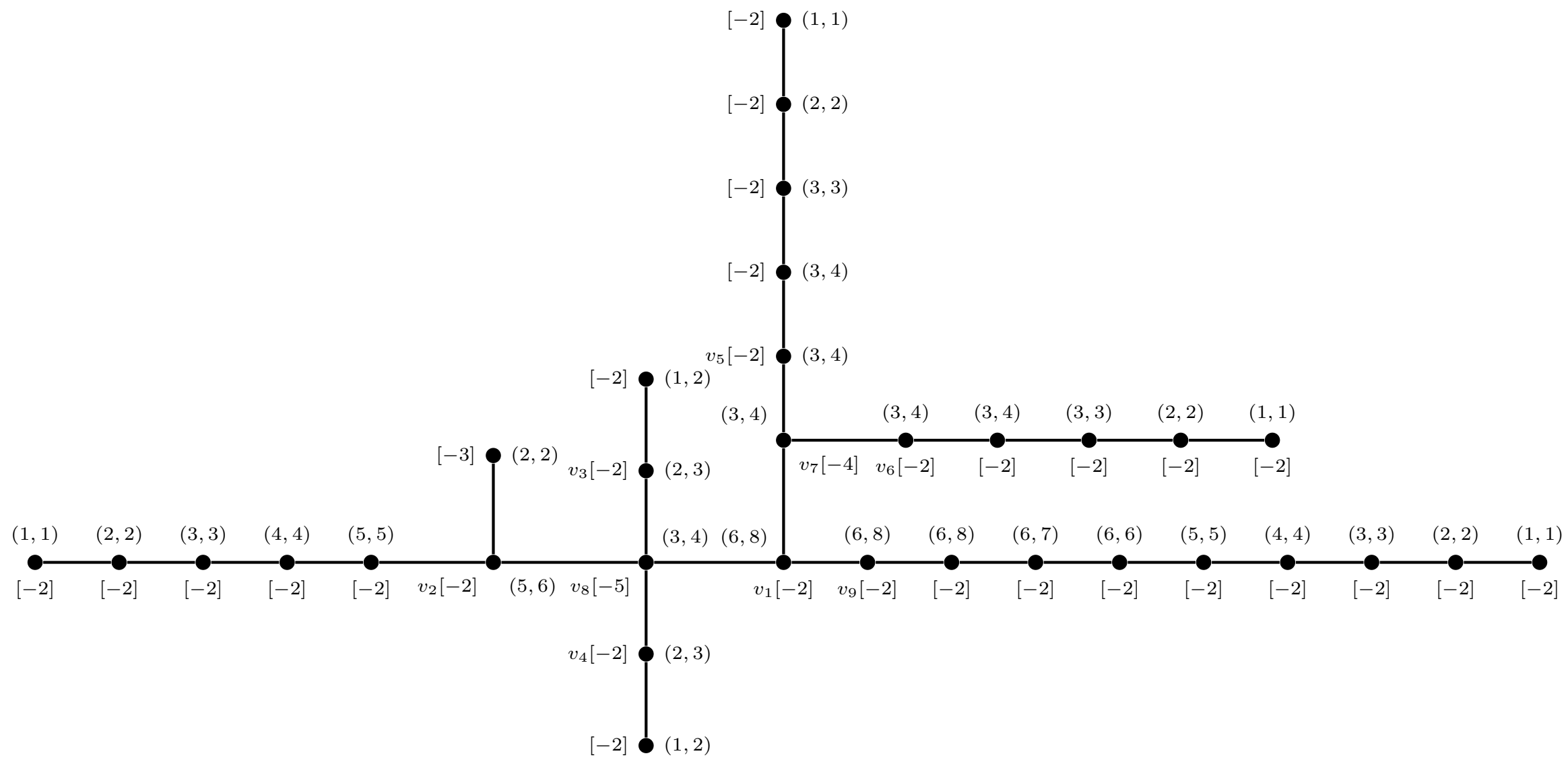
\begin{figure}[htbp]\label{grafmax}
  \centering
  \begin{tikzpicture}[scale=0.75, >=stealth, every node/.style={transform shape}]


    \node[circle, fill=black, inner sep=2pt, label=below:{\scriptsize $[-2]$}, label=above:{\scriptsize $(1,1)$}] (v1) at (-8.0, 0.0) {};
    \node[circle, fill=black, inner sep=2pt, label=below:{\scriptsize $[-2]$}, label=above:{\scriptsize $(2,2)$}] (v2) at (-6.9, 0.0) {};
    \node[circle, fill=black, inner sep=2pt, label=below:{\scriptsize $[-2]$}, label=above:{\scriptsize $(3,3)$}] (v3) at (-5.8, 0.0) {};
    \node[circle, fill=black, inner sep=2pt, label=below:{\scriptsize $[-2]$}, label=above:{\scriptsize $(4,4)$}] (v4) at (-4.7, 0.0) {};
    \node[circle, fill=black, inner sep=2pt, label=below:{\scriptsize $[-2]$}, label=above:{\scriptsize $(5,5)$}] (v5) at (-3.6, 0.0) {};

    \node[circle, fill=black, inner sep=2pt, label=below left:{\scriptsize $v_{2}[-2]$}, label=below right:{\scriptsize $(5,6)$}] (v6) at (-2.0, 0.0) {};

    \node[circle, fill=black, inner sep=2pt, label=left:{\scriptsize $[-3]$}, label=right:{\scriptsize $(2,2)$}] (v7) at (-2.0, 1.4) {};

    \node[circle, fill=black, inner sep=2pt, label=below left:{\scriptsize $v_{8}[-5]$}, label=above right:{\scriptsize $(3,4)$}] (v8) at (0.0, 0.0) {};

    \node[circle, fill=black, inner sep=2pt, label=left:{\scriptsize $v_3[-2]$}, label=right:{\scriptsize $(2,3)$}]   (v9)  at (0.0, 1.2) {};
    \node[circle, fill=black, inner sep=2pt, label=left:{\scriptsize $[-2]$}, label=right:{\scriptsize $(1,2)$}] (v10) at (0.0, 2.4) {};

    \node[circle, fill=black, inner sep=2pt, label=left:{\scriptsize $v_4[-2]$}, label=right:{\scriptsize $(2,3)$}] (v11) at (0.0, -1.2) {};
    \node[circle, fill=black, inner sep=2pt, label=left:{\scriptsize $[-2]$}, label=right:{\scriptsize $(1,2)$}] (v12) at (0.0, -2.4) {};

    \node[circle, fill=black, inner sep=2pt, label=below:{\scriptsize $v_{1}[-2]$}, label=above left:{\scriptsize $(6,8)$}] (v13) at (1.8, 0.0) {};

    \node[circle, fill=black, inner sep=2pt, label=below:{\scriptsize $v_{9}[-2]$}, label=above:{\scriptsize $(6,8)$}] (v14) at (2.9, 0.0) {};
    \node[circle, fill=black, inner sep=2pt, label=below:{\scriptsize $[-2]$}, label=above:{\scriptsize $(6,8)$}] (v15) at (4.0, 0.0) {};
    \node[circle, fill=black, inner sep=2pt, label=below:{\scriptsize $[-2]$}, label=above:{\scriptsize $(6,7)$}] (v16) at (5.1, 0.0) {};
    \node[circle, fill=black, inner sep=2pt, label=below:{\scriptsize $[-2]$}, label=above:{\scriptsize $(6,6)$}] (v17) at (6.2, 0.0) {};
    \node[circle, fill=black, inner sep=2pt, label=below:{\scriptsize $[-2]$}, label=above:{\scriptsize $(5,5)$}] (v18) at (7.3, 0.0) {};
    \node[circle, fill=black, inner sep=2pt, label=below:{\scriptsize $[-2]$}, label=above:{\scriptsize $(4,4)$}] (v19) at (8.4, 0.0) {};
    \node[circle, fill=black, inner sep=2pt, label=below:{\scriptsize $[-2]$}, label=above:{\scriptsize $(3,3)$}] (v20) at (9.5, 0.0) {};
    \node[circle, fill=black, inner sep=2pt, label=below:{\scriptsize $[-2]$}, label=above:{\scriptsize $(2,2)$}] (v21) at (10.6, 0.0) {};
    \node[circle, fill=black, inner sep=2pt, label=below:{\scriptsize $[-2]$}, label=above:{\scriptsize $(1,1)$}] (v22) at (11.7, 0.0) {};

    \node[circle, fill=black, inner sep=2pt, label=below right:{\scriptsize $v_{7}[-4]$}, label=above left:{\scriptsize $(3,4)$}] (v23) at (1.8, 1.6) {};

    \node[circle, fill=black, inner sep=2pt, label=left:{\scriptsize $v_{5}[-2]$}, label=right:{\scriptsize $(3,4)$}] (v24) at (1.8, 2.7) {};
    \node[circle, fill=black, inner sep=2pt, label=left:{\scriptsize $[-2]$}, label=right:{\scriptsize $(3,4)$}] (v25) at (1.8, 3.8) {};
    \node[circle, fill=black, inner sep=2pt, label=left:{\scriptsize $[-2]$}, label=right:{\scriptsize $(3,3)$}] (v26) at (1.8, 4.9) {};
    \node[circle, fill=black, inner sep=2pt, label=left:{\scriptsize $[-2]$}, label=right:{\scriptsize $(2,2)$}] (v27) at (1.8, 6.0) {};
    \node[circle, fill=black, inner sep=2pt, label=left:{\scriptsize $[-2]$}, label=right:{\scriptsize $(1,1)$}] (v28) at (1.8, 7.1) {};

    \node[circle, fill=black, inner sep=2pt, label=below:{\scriptsize $v_{6}[-2]$}, label=above:{\scriptsize $(3,4)$}] (v29) at (3.4, 1.6) {};
    \node[circle, fill=black, inner sep=2pt, label=below:{\scriptsize $[-2]$}, label=above:{\scriptsize $(3,4)$}] (v30) at (4.6, 1.6) {};
    \node[circle, fill=black, inner sep=2pt, label=below:{\scriptsize $[-2]$}, label=above:{\scriptsize $(3,3)$}] (v31) at (5.8, 1.6) {};
    \node[circle, fill=black, inner sep=2pt, label=below:{\scriptsize $[-2]$}, label=above:{\scriptsize $(2,2)$}] (v32) at (7.0, 1.6) {};
    \node[circle, fill=black, inner sep=2pt, label=below:{\scriptsize $[-2]$}, label=above:{\scriptsize $(1,1)$}] (v33) at (8.2, 1.6) {};


    \draw[thick] (v1) -- (v2);
    \draw[thick] (v2) -- (v3);
    \draw[thick] (v3) -- (v4);
    \draw[thick] (v4) -- (v5);
    \draw[thick] (v5) -- (v6);
    \draw[thick] (v6) -- (v8);
    \draw[thick] (v6) -- (v7);

    \draw[thick] (v8) -- (v9);
    \draw[thick] (v9) -- (v10);
    \draw[thick] (v8) -- (v11);
    \draw[thick] (v11) -- (v12);

    \draw[thick] (v8) -- (v13);
    \draw[thick] (v13) -- (v14);
    \draw[thick] (v14) -- (v15);
    \draw[thick] (v15) -- (v16);
    \draw[thick] (v16) -- (v17);
    \draw[thick] (v17) -- (v18);
    \draw[thick] (v18) -- (v19);
    \draw[thick] (v19) -- (v20);
    \draw[thick] (v20) -- (v21);
    \draw[thick] (v21) -- (v22);

    \draw[thick] (v13) -- (v23);

    \draw[thick] (v23) -- (v24);
    \draw[thick] (v24) -- (v25);
    \draw[thick] (v25) -- (v26);
    \draw[thick] (v26) -- (v27);
    \draw[thick] (v27) -- (v28);

    \draw[thick] (v23) -- (v29);
    \draw[thick] (v29) -- (v30);
    \draw[thick] (v30) -- (v31);
    \draw[thick] (v31) -- (v32);
    \draw[thick] (v32) -- (v33);

  \end{tikzpicture}
  \caption{The multiplicity system of $Z_{min}$  and ${\rm div}_E\pi^*y$ in $(\, , \, )$}
  \label{fig:resolution_graph_zmax_perfect}
\end{figure}

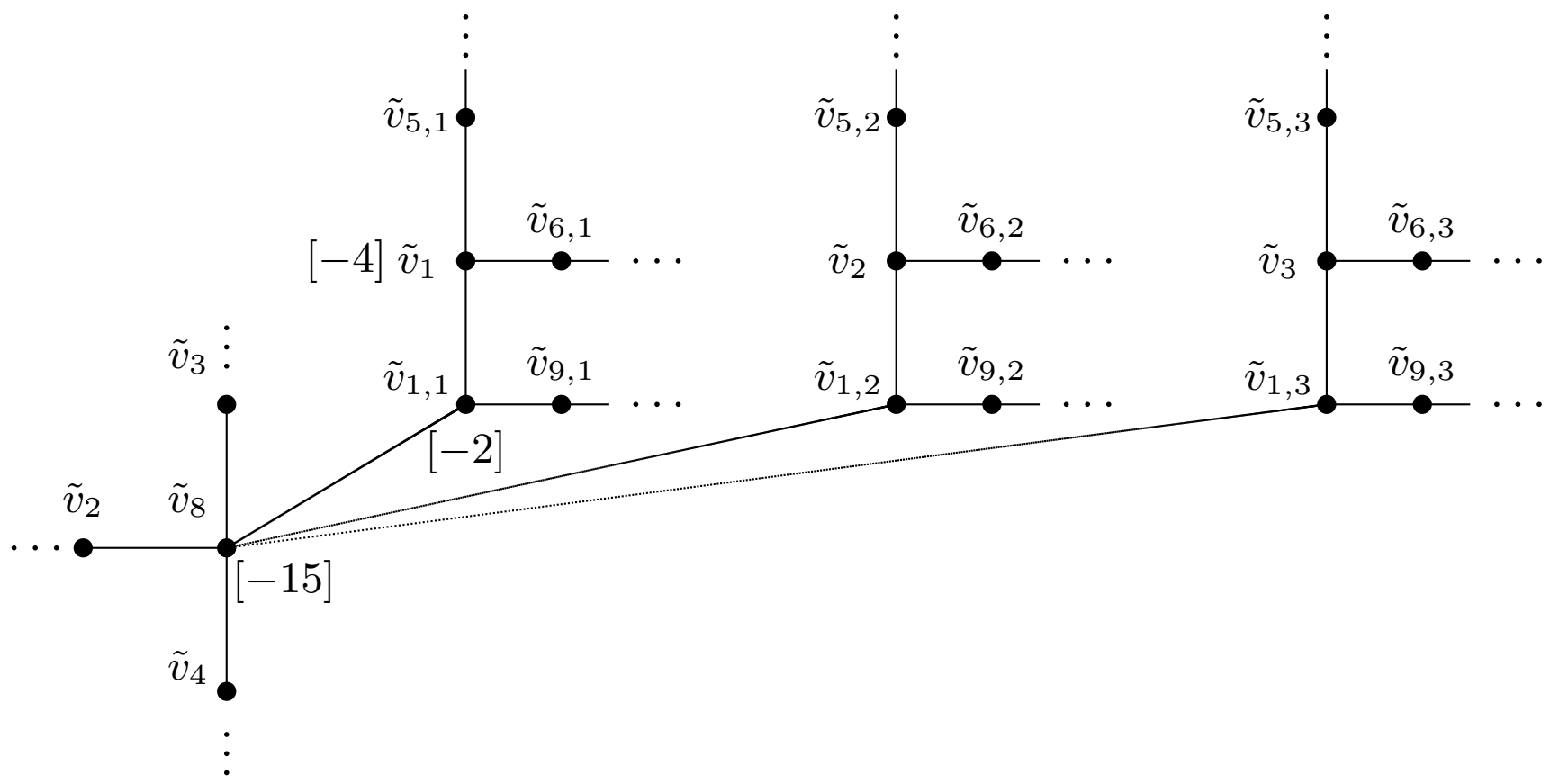
\begin{figure}[htbp]
 \begin{picture}(400,150)(0,0)


    \put(0,30){\circle*{4}}
    \put(30,0){\circle*{4}}
    \put(30,30){\circle*{4}}
    \put(30,60){\circle*{4}}

  \put(80,60){\circle*{4}}
    \put(80,90){\circle*{4}}
    \put(80,120){\circle*{4}}
  \put(100,60){\circle*{4}}
    \put(100,90){\circle*{4}}

 \put(170,60){\circle*{4}}
    \put(170,90){\circle*{4}}
    \put(170,120){\circle*{4}}
  \put(190,60){\circle*{4}}
    \put(190,90){\circle*{4}}

 \put(260,60){\circle*{4}}
    \put(260,90){\circle*{4}}
    \put(260,120){\circle*{4}}
  \put(280,60){\circle*{4}}
    \put(280,90){\circle*{4}}

    \put(0,30){\line(1,0){30}}
    \put(30,0){\line(0,1){60}}

\qbezier(30,30)(80,60)(80,60)
  \put(80,60){\line(1,0){30}}
  \put(80,90){\line(1,0){30}}
 \put(80,60){\line(0,1){70}}

\qbezier(30,30)(170,60)(170,60)
  \put(170,60){\line(1,0){30}}
  \put(170,90){\line(1,0){30}}
 \put(170,60){\line(0,1){70}}

\qbezier(30,30)(260,60)(260,60)
  \put(260,60){\line(1,0){30}}
  \put(260,90){\line(1,0){30}}
 \put(260,60){\line(0,1){70}}

    \put(-10,30){\makebox(0,0){$\ldots$}}
   \put(30,75){\makebox(0,0){$\vdots$}}
 \put(30,-10){\makebox(0,0){$\vdots$}}

   \put(80,140){\makebox(0,0){$\vdots$}}
 \put(170,140){\makebox(0,0){$\vdots$}}
 \put(260,140){\makebox(0,0){$\vdots$}}

   \put(120,60){\makebox(0,0){$\ldots$}}
   \put(120,90){\makebox(0,0){$\ldots$}}

   \put(210,60){\makebox(0,0){$\ldots$}}
   \put(210,90){\makebox(0,0){$\ldots$}}

   \put(300,60){\makebox(0,0){$\ldots$}}
   \put(300,90){\makebox(0,0){$\ldots$}}

    \put(0,40){\makebox(0,0){\small{$\tilde{v}_2$}}}
   \put(22,40){\makebox(0,0){\small{$\tilde{v}_8$}}}
   \put(22,70){\makebox(0,0){\small{$\tilde{v}_3$}}}
   \put(22,5){\makebox(0,0){\small{$\tilde{v}_4$}}}

   \put(70,65){\makebox(0,0){\small{$\tilde{v}_{1,1}$}}}
   \put(70,90){\makebox(0,0){\small{$\tilde{v}_1$}}}
   \put(70,120){\makebox(0,0){\small{$\tilde{v}_{5,1}$}}}
   \put(100,68){\makebox(0,0){\small{$\tilde{v}_{9,1}$}}}
  \put(100,98){\makebox(0,0){\small{$\tilde{v}_{6,1}$}}}

   \put(160,65){\makebox(0,0){\small{$\tilde{v}_{1,2}$}}}
   \put(160,90){\makebox(0,0){\small{$\tilde{v}_2$}}}
   \put(160,120){\makebox(0,0){\small{$\tilde{v}_{5,2}$}}}
   \put(190,68){\makebox(0,0){\small{$\tilde{v}_{9,2}$}}}
  \put(190,98){\makebox(0,0){\small{$\tilde{v}_{6,2}$}}}

   \put(250,65){\makebox(0,0){\small{$\tilde{v}_{1,3}$}}}
   \put(250,90){\makebox(0,0){\small{$\tilde{v}_3$}}}
   \put(250,120){\makebox(0,0){\small{$\tilde{v}_{5,3}$}}}
   \put(280,68){\makebox(0,0){\small{$\tilde{v}_{9,3}$}}}
  \put(280,98){\makebox(0,0){\small{$\tilde{v}_{6,3}$}}}

\put(42,23){\makebox(0,0){\small{$[-15]$}}}
\put(80,50){\makebox(0,0){\small{$[-2]$}}}
  \put(55,90){\makebox(0,0){\small{$[-4]$}}}

\end{picture}
  \caption{Part of the covering graph}
\end{figure}

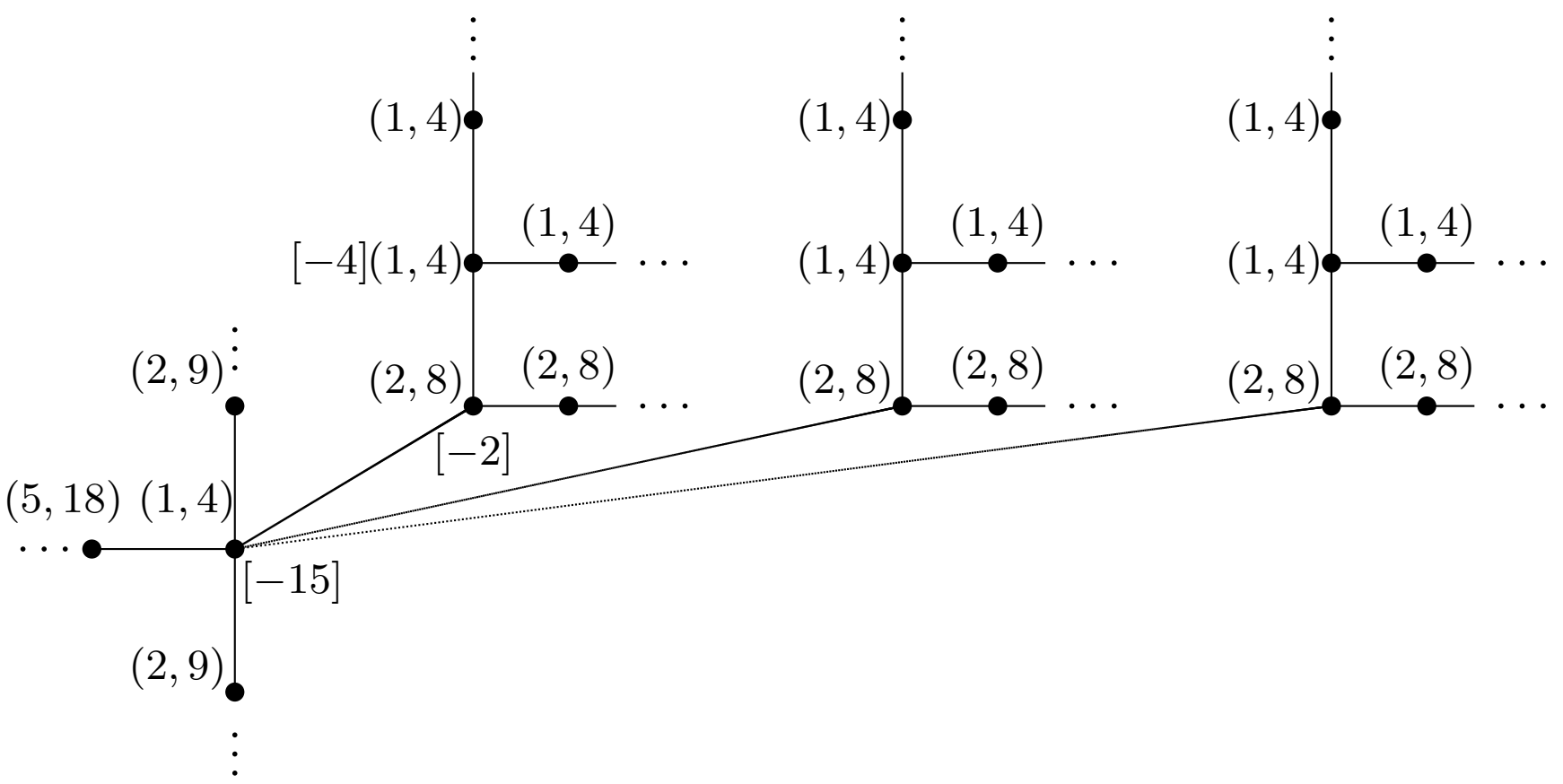
\begin{figure}[htbp]
 \begin{picture}(400,150)(0,0)


    \put(0,30){\circle*{4}}
    \put(30,0){\circle*{4}}
    \put(30,30){\circle*{4}}
    \put(30,60){\circle*{4}}

  \put(80,60){\circle*{4}}
    \put(80,90){\circle*{4}}
    \put(80,120){\circle*{4}}
  \put(100,60){\circle*{4}}
    \put(100,90){\circle*{4}}

 \put(170,60){\circle*{4}}
    \put(170,90){\circle*{4}}
    \put(170,120){\circle*{4}}
  \put(190,60){\circle*{4}}
    \put(190,90){\circle*{4}}

 \put(260,60){\circle*{4}}
    \put(260,90){\circle*{4}}
    \put(260,120){\circle*{4}}
  \put(280,60){\circle*{4}}
    \put(280,90){\circle*{4}}

    \put(0,30){\line(1,0){30}}
    \put(30,0){\line(0,1){60}}

\qbezier(30,30)(80,60)(80,60)
  \put(80,60){\line(1,0){30}}
  \put(80,90){\line(1,0){30}}
 \put(80,60){\line(0,1){70}}

\qbezier(30,30)(170,60)(170,60)
  \put(170,60){\line(1,0){30}}
  \put(170,90){\line(1,0){30}}
 \put(170,60){\line(0,1){70}}

\qbezier(30,30)(260,60)(260,60)
  \put(260,60){\line(1,0){30}}
  \put(260,90){\line(1,0){30}}
 \put(260,60){\line(0,1){70}}

    \put(-10,30){\makebox(0,0){$\ldots$}}
   \put(30,75){\makebox(0,0){$\vdots$}}
 \put(30,-10){\makebox(0,0){$\vdots$}}

   \put(80,140){\makebox(0,0){$\vdots$}}
 \put(170,140){\makebox(0,0){$\vdots$}}
 \put(260,140){\makebox(0,0){$\vdots$}}

   \put(120,60){\makebox(0,0){$\ldots$}}
   \put(120,90){\makebox(0,0){$\ldots$}}

   \put(210,60){\makebox(0,0){$\ldots$}}
   \put(210,90){\makebox(0,0){$\ldots$}}

   \put(300,60){\makebox(0,0){$\ldots$}}
   \put(300,90){\makebox(0,0){$\ldots$}}

    \put(-6,40){\makebox(0,0){\small{$(5,18)$}}}
   \put(20,40){\makebox(0,0){\small{$(1,4)$}}} \put(42,23){\makebox(0,0){\small{$[-15]$}}}
   \put(18,67){\makebox(0,0){\small{$(2,9)$}}}
   \put(18,5){\makebox(0,0){\small{$(2,9)$}}}

   \put(68,65){\makebox(0,0){\small{$(2,8)$}}} \put(80,50){\makebox(0,0){\small{$[-2]$}}}
   \put(68,90){\makebox(0,0){\small{$(1,4)$}}} \put(50,90){\makebox(0,0){\small{$[-4]$}}}
   \put(68,120){\makebox(0,0){\small{$(1,4)$}}}
   \put(100,68){\makebox(0,0){\small{$(2,8)$}}}
  \put(100,98){\makebox(0,0){\small{$(1,4)$}}}

   \put(158,65){\makebox(0,0){\small{$(2,8)$}}}
   \put(158,90){\makebox(0,0){\small{$(1,4)$}}}
   \put(158,120){\makebox(0,0){\small{$(1,4)$}}}
   \put(190,68){\makebox(0,0){\small{$(2,8)$}}}
  \put(190,98){\makebox(0,0){\small{$(1,4)$}}}

   \put(248,65){\makebox(0,0){\small{$(2,8)$}}}
   \put(248,90){\makebox(0,0){\small{$(1,4)$}}}
   \put(248,120){\makebox(0,0){\small{$(1,4)$}}}
   \put(280,68){\makebox(0,0){\small{$(2,8)$}}}
  \put(280,98){\makebox(0,0){\small{$(1,4)$}}}

\end{picture}
\caption{Part of the   multiplicity system of $z$ and $f^*_3y$  on $\tX^{(3)}$ }
\end{figure}

\end{document}